\documentclass[11pt,a4paper]{article}

\usepackage[T1]{fontenc}
\usepackage[utf8]{inputenc}
\usepackage[a4paper,margin=1in]{geometry}
\usepackage{amsmath,amssymb,amsfonts,mathrsfs}
\usepackage{amsthm}
\usepackage{enumerate}
\usepackage{xcolor}
\usepackage[authoryear,longnamesfirst]{natbib}
\usepackage[colorlinks=true,linkcolor=blue,citecolor=blue,urlcolor=blue]{hyperref}

\definecolor{revisionred}{rgb}{0,0,0}
\definecolor{REVISIONRED}{rgb}{0,0,0}
\definecolor{RED}{rgb}{0,0,0}
\definecolor{revisionblue}{rgb}{0,0,0}
\definecolor{reviewpurple}{rgb}{0,0,0}
\definecolor{originalgray}{rgb}{0,0,0}

\theoremstyle{plain}

\newtheorem{theorem}{Theorem}[section]

\newtheorem{tm}[theorem]{Theorem}
\newtheorem{prop}[theorem]{Proposition}
\newtheorem{lm}[theorem]{Lemma}

\theoremstyle{definition}

\newtheorem{ap}[theorem]{Assumption}

\theoremstyle{remark}

\newtheorem{rk}[theorem]{Remark}

\numberwithin{equation}{section}
\newcommand{\E}{\mathbb E}

\newcommand{\N}{\mathbb N}

\newcommand{\HH}{\mathbb H}

\newcommand{\norm}[1]{\left\|#1\right\|}
\newcommand{\pie}{\pi_{N_\epsilon}}

\newcommand{\ip}[2]{\left\langle #1,#2\right\rangle}
\newcommand{\<}{\langle}
\renewcommand{\>}{\rangle}
\newcommand{\rev}[1]{#1}
\newcommand{\revblue}[1]{#1}
\newcommand{\revred}[1]{#1}
\newcommand{\questionred}[1]{}
\newcommand{\revpurple}[1]{#1}
\newcommand{\questionpurple}[1]{}

\title{Polynomial weak approximation for stochastic reaction-diffusion equations near the sharp interface limit\thanks{The research is partially supported by NSFC grant 12522119, NSFC grant 12301526, MOST National Key \& Program No. 2024YFA1015900, the Hong Kong Research Grant Council GRF grant 15302823, NSFC/RGC Joint Research Scheme N\_PolyU5141/24, and the CAS AMSS-PolyU Joint Laboratory of Applied Mathematics. This work is supported by the National Natural Science Foundation of China (No. 12471386, No. 12461160278).}}

\author{
Jianbo Cui\\
Department of Applied Mathematics, The Hong Kong Polytechnic University\\
Hung Hom, Kowloon, Hong Kong\\
\texttt{jianbo.cui@polyu.edu.hk}
\and
Liying Sun\\
LSEC, ICMSEC, Academy of Mathematics and Systems Science\\
Chinese Academy of Sciences\\
Beijing, China\\
\texttt{liyingsun@lsec.cc.ac.cn}
}

\date{}

\begin{document}
\allowdisplaybreaks
\maketitle

\begin{abstract}
We study weak approximation for stochastic reaction-diffusion equations in the sharp-interface regime, where the diffuse interface thickness \(\epsilon\) is small and the dependence of numerical constants on \(\epsilon^{-1}\) is a central issue. For an additive-noise stochastic Allen--Cahn type equation, direct stability arguments typically produce weak error bounds with constants growing exponentially in \(\epsilon^{-1}\). Such estimates do not capture the polynomial stability expected near the sharp interface limit.

We prove polynomial-in-\(\epsilon^{-1}\) weak error bounds for an accelerated splitting exponential Euler approximation. The proof combines time-uniform moment estimates, regularity estimates for the exact and numerical dynamics, and time-independent derivative estimates for the Kolmogorov semigroup. The key point is that the averaged regularizing effect of the noise, expressed through asymptotic strong Feller or strong Feller estimates, replaces the deterministic spectral estimate for the linearized Allen--Cahn operator. The result gives a weak approximation theory whose constants depend polynomially on \(\epsilon^{-1}\) and explicitly on the covariance regularity and non-degeneracy parameters.
\end{abstract}



\medskip

\section{Introduction}
\revblue{The purpose of this paper is to develop a weak approximation theory for a class of stochastic reaction-diffusion equations in a regime where the usual stability constants are not informative.  More precisely, we consider equations close to the sharp interface limit and ask whether the constants in numerical error estimates can be controlled by powers of the inverse interface thickness \(\epsilon^{-1}\), rather than by exponential functions of \(\epsilon^{-1}\).  This question is both analytic and numerical: it concerns the regularity of infinite-dimensional Markov semigroups generated by monotone polynomial drifts, and it determines whether a numerical error bound remains meaningful for small \(\epsilon\in (0,1)\).}

\revblue{The model problem is}
\begin{align}\label{sac}
	d u(t)=A u(t)dt-\frac 1 \epsilon f(u(t)) dt+dW(t).
\end{align}
\revblue{Here \(A\) is the Laplacian on a bounded strongly regular domain \(\mathcal D\subset\mathbb R^d\), \(d\le3\), with homogeneous Dirichlet or Neumann boundary condition, \(W\) is an \(\mathbb H\)-valued \(Q\)-Wiener process on \(\mathbb H=L^2(\mathcal D)\), and \(f\) is the Nemytskii operator associated with an odd-degree polynomial.  The Allen--Cahn nonlinearity \(f(\xi)=\xi^3-\xi\) is the guiding example.}

\revblue{Strong and weak approximations of semilinear parabolic SPDEs with globally Lipschitz or sufficiently regular coefficients are now well developed; see, for instance, \cite{Deb11,JKW11,KLL12,KLL13,BL13,AKL16,BD18} for representative weak approximation results and \cite{Gyo99,DP09,LPS14,Yan05} for related strong approximation results.  For superlinear drifts, stable discretizations and convergence estimates require additional ideas; see, among others, \cite{BJ16,BG18,BG18B,BCH18,CH18,CH20,CHS21,CHS21b,FLZ17,LQ20,LQ21,MP17,Wang20,QW20}.  The present work addresses a complementary issue that is not captured by a fixed-\(\epsilon\) convergence theory: the explicit dependence of weak error constants on the singular parameter \(\epsilon^{-1}\).}

\revblue{This dependence is delicate already at a formal level.  For a temporal approximation of the stochastic Allen--Cahn equation, the one-sided Lipschitz constant of \(\epsilon^{-1}(-\xi^3+\xi)\) is of order \(\epsilon^{-1}\).  A direct Gronwall argument therefore leads to an estimate of the form}
\begin{align*}
	\left\|u(t_k)-u_k\right\|^2_{L^2(\Omega;\mathbb H)} \le \frac C{\epsilon} \sum_{j=0}^{k-1}\left\|u(t_{j})-u_{j}\right\|_{L^2(\Omega;\mathbb H)}^2\tau +\textbf{Err}
\end{align*}
\revblue{for \(T=N\tau\), \(t_j=j\tau\), $k\le N$, and an error term \(\textbf{Err}\to0\) as \(\tau\to0\).  Consequently, \(\|u(T)-u_N\|_{L^2(\Omega;\mathbb H)}\) is controlled by a factor of order \(\exp(CT/\epsilon)\).  The same exponential dependence is inherited by weak errors if one bounds}
\begin{align*}
	\Big|\E \Big[\phi(u(T))-\phi(u_N)\Big]\Big|\le C|D\phi|_{\infty}\|u(t_k)-u_k\|_{L^2(\Omega;\mathbb H)}
\end{align*}
for \(\phi\in\mathcal C_b^1(\mathbb H)\).  The central difficulty is therefore to avoid importing this exponential factor into either strong or weak numerical error estimates.

\revblue{In deterministic Allen--Cahn theory, polynomial-in-\(\epsilon^{-1}\) estimates are often obtained from spectral estimates for the linearized Allen--Cahn operator; see, for example, \cite{Chen1994,FP03}.  However, such estimates are not directly available for stochastic equations with additive noise.  A recent result \cite{ABNP21} obtains polynomial dependence for an implicit approximation in a small-noise regime by combining deterministic spectral information with perturbative arguments.  Here we follow a different route: we use the smoothing of the infinite-dimensional Markov semigroup itself and do not require the noise intensity to vanish with \(\epsilon\).}

\revblue{The main analytic contributions of this paper are the following.}
\begin{enumerate}
	\item[(a)] \revred{We establish time-uniform moment and regularity estimates for the original equation \eqref{sac} and for the accelerated splitting approximation \eqref{spl-sch}.}
	\item[(b)] \revblue{We prove time-independent first- and second-order derivative estimates for the Kolmogorov semigroup \(P_t\phi(x)=\E[\phi(u(t,x))]\).  The constants are polynomial in \(\epsilon^{-1}\) and explicit in the covariance regularity quantity \(\mathfrak q_\beta\), the low-mode non-degeneracy quantity \(\gamma(\epsilon)\), and, under full non-degeneracy, the inverse covariance norm \(\mathfrak q_-\).}
	
	\item[(c)] 
	\revblue{We combine these estimates with a weak error expansion for an accelerated splitting exponential Euler method and obtain weak convergence with constants depending at most polynomially on \(\epsilon^{-1}\).  The proof is written for a temporal discretization, but the decomposition separates the semigroup regularity from the particular scheme and indicates how related spatial or full discretizations may be treated once analogous stability estimates are available.}   
	
\end{enumerate}
\section{Main results}
\label{sec-main}

\revblue{This section contains the assumptions, the notation for the covariance-dependent constants, the Kolmogorov regularity theorem, and the weak convergence theorem.  The order is chosen so that the analytic input for the weak error estimate is visible before the numerical scheme is introduced.}

We first state the notation and assumptions that are needed for the main results.
\revblue{The parameter \(\lambda\ge 0\) is fixed throughout the paper and is not allowed to depend on \(\epsilon^{-1}\).} Throughout this paper, we formulate the
time-step restriction simply as
\[
\qquad 0<\tau\lesssim \epsilon .
\]
All constants below are independent of \(N\), \(\tau\), and \(\frac 1\epsilon\), unless
their dependence is explicitly indicated.

\subsection{Preliminaries}

\revblue{We use the following notation throughout the paper:}
\[
\revblue{\mathbb H=L^2(\mathcal D),\qquad \|\cdot\|=\|\cdot\|_{\mathbb H},\qquad E=C(\overline{\mathcal D}).}
\]
\revblue{Unless another norm is explicitly displayed, \(\|\cdot\|\) denotes the \(\mathbb H\)-norm.  Given two separable Hilbert spaces \(\mathcal H\) and \(\widetilde H\), let \(\mathcal L(\mathcal H,\widetilde H)\) be the space of bounded linear operators from \(\mathcal H\) to \(\widetilde H\).  We denote by \(HS(\mathcal H,\widetilde H)\) the Hilbert--Schmidt space equipped with}
\[
\revblue{\|\cdot\|_{HS(\mathcal H,\widetilde H)}
:=\left(\sum_{j\in\mathbb N^+}\|\cdot (g_j)\|_{\widetilde H}^2\right)^{1/2},}
\]
\revblue{where \(\{g_j\}_{j\in\mathbb N^+}\) is any orthonormal basis of \(\mathcal H\).  When \(\mathcal H=\widetilde H=\mathbb H\), we write \(HS(\mathbb H)\).}

Throughout the paper, we assume the following conditions unless otherwise stated.

\begin{ap}\label{ass-main}
\begin{enumerate}[(A1)]
\item[(i)] \rev{The operator \(A\) is the Laplacian on \(\mathcal D\) with homogeneous Dirichlet or Neumann boundary condition.  The covariance operator \(Q\) is bounded, non-negative, self-adjoint, and commutes with \(A\).}
\item[(ii)] \rev{The drift is the Nemytskii operator of  a polynomial}
\[
\rev{f(\xi)=\sum_{j=0}^{2m+1}a_j\xi^j,\qquad f(0)=0,\qquad a_{2m+1}>0.}
\]
\end{enumerate}
\rev{
The condition (ii) ensures that \(-\epsilon^{-1}f\) is dissipative at infinity and locally Lipschitz. 
The Allen--Cahn-type choice \(f(\xi)=\xi^3-\xi\) is the canonical example.
Note that a bounded, smooth, Lipschitz perturbation may be added to the drift without changing the argument.  We do not include this extension in the paper.}
\end{ap}

  \revred{By condition (i), one can choose  \(\lambda>0\)  such that \(-A+\lambda I\) is a strictly positive operator generating an analytic semigroup on \(L^p:=L^p(\mathcal D)\), \(2\le p<+\infty\) (see e.g. \cite[Chapter 6]{Cer01}). Under the homogeneous Neumann boundary condition, \(-A+\lambda I\) generates an analytic and strongly continuous semigroup on \(E=\mathcal C(\overline{\mathcal D})\); under the Dirichlet boundary condition, it generates an analytic and strongly continuous semigroup on the subspace of \(E\) consisting of functions satisfying the zero boundary condition.}

 We will frequently use the following smoothing effect of $S_{\lambda}(t)=e^{(A-\lambda)t}$, that is,   for any $\kappa\ge 0,$
and $c_0\in [0,\lambda)$ (when $\lambda=0,$ we set $c_0=0$), there  exists a constant $C_{\kappa}>0$ such that  \begin{align}\label{smooth0}
	{\left\|S_{\lambda}(t)(-A+\lambda I)^{\kappa }\right\|_{\mathcal L(\mathbb H)}
	\le C_{\kappa}e^{-c_0t}t^{-\kappa },}
\end{align}
and the continuity estimate (see e.g. \cite[Appendix B]{Kru14a}), i.e.,  for any $\kappa_1\in [0,1],$
\begin{align}\label{smooth1}
	&{\left\|(-A+\lambda I)^{-\kappa_1}(S_{\lambda}(t)-I)\right\|_{\mathcal L(\mathbb H)} \le C_{\kappa_1} t^{\kappa_1}, \; t\ge 0,}\\\label{smooth2}
	&{\Big\|(-A+\lambda I)^{\kappa_1}\int_{t_1}^{t_2} S_{\lambda}(t_2-\sigma) v\,d\sigma\Big\|
	\le C_{\kappa_1}(t_2-t_1)^{1-\kappa_1}\|v\|, \; t_2\ge t_1\ge 0.}
\end{align}
Moreover, the super-contractility property (see e.g. \cite[Appendix B]{PZ07}) holds, that is,  for any $c_0\in [0,\lambda)$ and $1\le p\le q<\infty,$  there exist positive constants $C_{p,q}$ and $C_p$ such that  
\begin{align}\label{super0}
	&\left\|S_{\lambda}(t)v\right\|_{L^q}\le C_{p,q}e^{-c_0 t}t^{-\frac d2\left(\frac 1p-\frac 1q\right)}\|v\|_{L^p}, \; \text{for} \; v\in L^p,\\\label{super1}
	&\left\|S_{\lambda}(t)v\right\|_E\le C_{p}e^{-c_0 t}t^{-\frac d{2p}}\|v\|_{L^p}, \; \text{for} \;v\in L^p.
\end{align}
Denote $\mathbb H^{\kappa}$ as  $D((-A+\lambda)^{\frac \kappa 2})$ with the norm 
\begin{align*}
	\left\|v\right\|_{\mathbb H^{\kappa}}=\Big(\sum_{j\in \mathbb N^+}(\lambda_j+\lambda)^{\kappa} \left\<v,f_j\right\>^2\Big)^{\frac 12} 
\end{align*}
for $\kappa \ge 0$. Here $(f_j,\lambda_j)_{j\in \mathbb N^+}$ is an orthonormal eigensystem of $A$.   It is known that the norm $\|\cdot\|_{\mathbb H^{\kappa}}$ is equivalent to the standard norm of the Sobolev space $H^{\kappa}=W^{\kappa,2}(\mathcal D)$ equipped with the corresponding boundary condition induced by the realization of $A$. 
For $\kappa<0$, one could also define $\mathbb H^{\kappa}$ by the set of the series sum $v=\sum_{j}v_je_j$ satisfying that the norm $\|v\|_{\mathbb H^{\kappa}}:=(\sum_{j\in \mathbb N^+}(\lambda_j+\lambda)^{\kappa}v_j^2)^{\frac 12}$ is finite.  It can be seen that $\mathbb H^0=\mathbb H.$

Below, we present several useful settings and notations on the covariance operator $Q$ of the driving Wiener process. 
For \(\beta\in[0,2)\), we write
\[
	\mathfrak q_\beta
	:=1+\revblue{\|Q\|_{\mathcal L(\mathbb H)}}
	+\left\|(-A+\lambda)^{\frac{\beta-1}{2}}Q^{\frac12}\right\|_{HS(\mathbb H)}.
\]
When \(Q^{-1/2}\) exists and is bounded, we also set
\[
	\mathfrak q_-:=1+\left\|Q^{-\frac12}\right\|_{\mathcal L(\mathbb H)}.
\]
For \(q\in\mathbb N\), let \(V(x)=1+\|x\|^q\) and define
\[
	\mathcal B_{V,q}(\mathbb H)
	:=
	\left\{\phi:\mathbb H\to\mathbb R \hbox{ Borel measurable}:
	\|\phi\|_{V,q}:=\sup_{x\in\mathbb H}
	\frac{|\phi(x)|}{1+\|x\|^q}<\infty
	\right\}.
\]
When \(q=0\), we write \(\mathcal B_b(\mathbb H)=\mathcal B_{V,0}(\mathbb H)\).

For \(\epsilon\in(0,1)\), set
\begin{equation}\label{eq:Neps}
    N_\epsilon
    :=\inf\left\{N\in\N:\lambda_{N+1}\ge \frac{1+L_f}{\epsilon}\right\},
    \qquad
    \gamma(\epsilon):=\sup_{1\le i\le N_\epsilon}q_i^{-1/2},
\end{equation}
where \(q_i\) are the eigenvalues of \(Q\) in the common eigenbasis of \(A\) and
\(Q\), and \(L_f>0\) is the one-sided Lipschitz constant of \(f\).

Assume that \(\gamma(\epsilon)<\infty\), which is a low-mode non-degeneracy assumption that will be used to establish the asymptotic strong Feller estimate. If $Q^{-\frac 12}$ exists, the additional condition \begin{align}\label{erg-noi}
	\left\|Q^{-\frac12}\right\|_{\mathcal L(\mathbb H)}<\infty ,
\end{align} will be used to derive the second-order regularity estimates of the Kolmogorov equation when $Q^{-\frac 12}$ exists. Note that the full non-degeneracy condition \eqref{erg-noi} is very restrictive when combined with
\(\|(-A+\lambda)^{(\beta-1)/2}Q^{1/2}\|_{HS(\mathbb H)}<\infty\). If \(q_i\) are bounded below by a positive constant, Weyl's law gives the compatibility condition \(\beta<1-d/2\). Thus it is empty for \(d\ge2\), and in \(d=1\) it requires \(\beta<1/2\). 

\revblue{Throughout the main results, a positive polynomial in displayed
variables means an ordinary finite polynomial with positive coefficients, whose
coefficients may depend on fixed model parameters such as \(f,m,\lambda,\beta\),
but not on \(N,\tau,\epsilon\), or the norms $\mathfrak q_\beta, \mathfrak q_-$.  In particular, the subindex of a polynomial denotes its degree with respect to the corresponding variable.}

\subsection{Main result I: regularity of the Kolmogorov semigroup}

We denote by $\mathcal C_b^k(\mathbb H)$, $k\in \mathbb N^+,$ the space of $k$ times continuously differentiable functionals from $\mathbb H$ to $\mathbb R$ whose derivatives up to order $k$ are bounded. Define
\begin{align*}
	&|D\phi|_\infty=\sup_{x\in \mathbb H} \|D\phi(x)\|_{\mathbb H}, |D^2\phi|_\infty=\sup_{x\in \mathbb H} \left\|D^2\phi(x)\right\|_{\mathcal L(\mathbb H)}.
\end{align*}
Here \(\revred{D^k\phi}\), \(k\le2\), denotes \revred{the \(k\)-th} derivative of $\phi\in \mathcal C_b^k(\mathbb H)$.

Let $\phi: \mathbb H \to \mathbb R$ be a function of class $\mathcal C_b^2(\mathbb H)$, and define
\[
X(t,x)=P_t\phi(x)=\mathbb E_x[\phi(u(t))]=\mathbb E[\phi(u(t,x))].
\]
Formally, \(X\) solves the following Kolmogorov equation \cite{Cer01}
\begin{align*} 
	\partial_t X(t,x)&=D X(t,x)\cdot \left[Ax-\frac 1\epsilon f(x)\right] +\frac 12 \sum_{j\in \mathbb N^+} D^2 X(t,x)\cdot \left(Q^{\frac 12}e_j,Q^{\frac 12}e_j\right).
\end{align*}
\iffalse
More precisely, define $X^M: [0,+\infty)\times H \to \mathbb R$ by $X^M(t,x)=\mathbb E_{x} [\phi(\mathcal P^M u(t))]$
with $u(0)=x\in H$, where $M$ is the parameter of the spectral Galerkin projection $\mathcal P^M$, i.e., $P^Mw=\sum_{i=1}^M \<w,e_i\>e_i$. Similar to \cite{}, one can verify that $X^M$ is the solution of  following regularized Kolmogorov equation
\begin{align*}
	\partial_t X^M(t,x)&=D X^M(t,x)\cdot[A\mathcal P^Mx+\mathcal P^M f(x)] \\
	&+\frac 12 \sum_{j\in \mathbb N^+} D^2 X^M(t,x)\cdot(\mathcal P_MQ^{\frac 12}e_j,\mathcal P_MQ^{\frac 12}e_j).
\end{align*}
\fi
 Notice that for any 
$h,h_1,h_2\in \mathbb H$, it holds that 
\begin{align}\label{rep-kol}
	DX(t,x)\cdot h&=\E_{x} \left[D\phi(u(t))\cdot \eta^h(t)\right],\\\label{rep-kol1}
	D^2X(t,x)\cdot (h_1,h_2)&=\E_{x} \left[D\phi(u(t))\cdot \zeta^{h_1,h_2}(t)\right]\\\nonumber
	&+\E_{x}\left[D^2\phi({u(t)})\cdot (\eta^{h_1}(t),\eta^{h_2}(t))\right].
\end{align}
Here $\eta^h$ and $\zeta^{h_1,h_2}$ satisfy the first order variational equation 
\begin{align}\label{var-1}
	\revred{\partial_t \eta ^h= A\eta ^h- \frac 1\epsilon f'(u(t))\cdot \eta^h,\qquad \eta^h(0)=h,}
\end{align}
and the second order variational equation 
\begin{align}\label{var-2}
	\revred{\partial_t \zeta^{h_1,h_2}}&\revred{=A\zeta ^{h_1,h_2}-\frac 1\epsilon f'(u(t))\cdot \zeta ^{h_1,h_2}
	-\frac 1\epsilon f''(u(t))\cdot \left(\eta^{h_1},\eta^{h_2}\right),}\\\nonumber
	\revred{\zeta^{h_1,h_2}(0)}&\revred{=0,}
\end{align}
respectively. 

\revblue{Our first main result is a time-uniform regularity estimate for the Kolmogorov semigroup.  This is the analytic input that prevents exponential dependence on \(\epsilon^{-1}\) in the weak error proof.  The estimate is stated directly for the original semigroup \(P_t\).  A
fully rigorous differentiability argument can be obtained by first proving the
same estimates for a regularized approximation, with bounds uniform in \(\delta\), and then passing to
the limit by the strong convergence of regularized approximations  \cite{CH18}. For simplicity, we omit these tedious details here.
}

\begin{tm}[Time-independent Kolmogorov regularity]\label{main-kol-reg}
	Let \(x\in \mathbb H^{\beta}\cap E\) and
	\(\|(-A+\lambda)^{(\beta-1)/2}Q^{1/2}\|_{HS(\mathbb H)}<\infty\) with
	\(\beta<2\). Assume
	\(\gamma(\epsilon)<\infty\).  Then, for any
	\(\alpha_1\in[0,1)\) and \(\phi\in\mathcal B_{V,q}(\mathbb H)\cap
	\mathcal C_b^1(\mathbb H)\), \(q\in\mathbb N\), there exist a positive polynomial
	\(P_{2m(2m+1)+q^2,2m}(\mathfrak q_\beta,\|x\|_E)\) and a constant $c_1\ge 0$ such that, for all \(t>0\) and
	\(h\in\mathbb H\),
	\begin{align}\label{time-ind1}
		|D P_t\phi(x)\cdot h|
		&\le P_{2m(2m+1)+q^2,2m}(\mathfrak q_\beta,\|x\|_E)(\|\phi\|_{V,q}+c_1\|D\phi\|_{\infty})
		a(\epsilon^{-1})
		\\\nonumber
		&\times \|h\|_{\mathbb H^{-2\alpha_1}}(1+t^{-\alpha_1}),
	\end{align}
	where \(a(\epsilon^{-1})=\gamma(\epsilon)\epsilon^{-5/2}\).
	
	If, in addition, \eqref{erg-noi} holds, then for any  
	\(\phi\in\mathcal B_{V,q}(\mathbb H), q\in \mathbb N\), 
	\eqref{time-ind1} holds with $c_1=0$ and $$a(\epsilon^{-1})=\mathfrak q_-\epsilon^{-5/2}.$$
	Furthermore, for any \(\phi\in\mathcal B_{V,q}(\mathbb H), q\in \mathbb N\),  and 
	\(\alpha_1,\alpha_2\in[0,1)\) with \(\alpha_1+\alpha_2<1\), 
	there exists a positive polynomial \( P_{(12m-1)(2m+1)+4m(8m-1)+q^2,10m-1}(\mathfrak q_\beta,\|x\|_E)\) such that for any $h_1,h_2\in \mathbb H,$
	\begin{align}\nonumber
	 |D^2P_t\phi(x)\cdot(h_1,h_2)|
	 &\le
	\|\phi\|_{V,q} \revblue{b(\epsilon^{-1})}\hat  P_{(12m-1)(2m+1)+4m(8m-1)+q^2,10m-1}(\mathfrak q_\beta,\|x\|_E)\\\label{time-ind2}
	 &\quad \|h_1\|_{\mathbb H^{-2\alpha_1}}\|h_2\|_{\mathbb H^{-2\alpha_2}} (1+t^{-\alpha_1-\alpha_2}),
	\end{align}
	where $b(\epsilon^{-1})= \mathfrak q_{-}\frac 1{\epsilon^{\frac {17}2}}$.
\end{tm}

\revblue{The condition \(\gamma(\epsilon)<\infty\) is a low-mode non-degeneracy condition.  We keep the \(\epsilon\)-dependence of \(\gamma(\epsilon)\) explicit because this quantity is part of the final polynomial dependence of the weak error constant (see Remark \ref{rk-gamma} for an example).  Under the stronger condition \eqref{erg-noi}, the Bismut--Elworthy--Li argument can be used, and the test function in \eqref{time-ind1}--\eqref{time-ind2} only needs polynomial growth and measurability.  The resulting second-order estimate \eqref{time-ind2}, with explicit dependence on \(\mathfrak q_\beta\), \(\mathfrak q_-\), and \(\epsilon^{-1}\), is one of the key new analytic inputs of the paper.}

\subsection{\texorpdfstring{Main result II: weak convergence with polynomial \(\epsilon^{-1}\)-dependence}{Main result II: weak convergence with polynomial epsilon inverse dependence}}

\revblue{We now introduce the temporal scheme used in the weak error theorem.} The accelerated splitting exponential Euler scheme \cite{BG18,BCH18} studied below is
\begin{align}\label{spl-sch}
 u_{n+1}=S_{\lambda}({\tau})\Phi_{\tau}^{\lambda}(u_n)
 +\int_{t_n}^{t_{n+1}} S_{\lambda}(t_{n+1}-s)dW(s),
 \qquad u_0=u(0),
\end{align}
where $\tau>0$ is the time step, \(t_n=n\tau\) for all \(n\in\mathbb N\). Here \(\Phi_t^\lambda\), $\lambda\ge 0,$ is the phase flow of the scalar equation
\[
    dX(t)=-\epsilon^{-1}f(X(t))\,dt+\lambda X(t)\,dt.
\]
\revred{The scheme \eqref{spl-sch} can be viewed as an exponential Euler scheme for the modified problem associated with \eqref{sac}, namely}
\begin{align}\label{reg-problem}
	\rev{dx^{\tau}(t)+(-A+\lambda)x^{\tau}(t) dt=\Psi_{\tau}^{\lambda}(x^{\tau}(t))dt+dW(t), \; x^{\tau}(0)=u(0).}
\end{align}
Here $\Psi_{\tau}^{\lambda}(x)=\frac {\Phi_{\tau}^{\lambda}(x)-x} {\tau}$. 
It has been shown in \cite{BCH18,BG18,CH18}, that the above regularizing problem \eqref{reg-problem} is a strong approximation of order $1$ for the original problem \eqref{sac} at any finite time when $\tau \to 0$. 

\revblue{For convenience, we separate two types of scalar flow estimates.  The first lemma contains local derivative and consistency bounds, while the second records dissipative bounds used for uniform moment estimates.}

\begin{lm}[Local flow estimates]\label{flow-local-est} Let $\lambda\ge 0$. For $t\ge0$ and $\xi\in \mathbb R$, the following estimates hold:
	\begin{align}\label{prop-phi}
		(\Phi_t^{\lambda})'(\xi)&\in [0, e^{c\left(\frac 1{\epsilon}+\lambda\right)t}], 
		\left|(\Phi_t^{\lambda})''(\xi)\right|\le C\left(\frac 1{\epsilon}+\lambda\right)e^{c\left(\frac 1{\epsilon}+\lambda\right)t }t\left(1+|\xi|^{2m-1}\right),\\\label{prop-psi}
		(\Psi_{t}^{\lambda})'(\xi)&\le  C\left(\frac 1{\epsilon}+\lambda\right)e^{c\left(\frac 1{\epsilon}+\lambda\right)t}, \; \left|(\Psi_{t}^{\lambda})'(\xi)\right| \le  C\left(\frac 1{\epsilon}+\lambda\right)e^{c\left(\frac 1{\epsilon}+\lambda\right)t}\left(1+|\xi|^{2m}\right),\\\label{prop-psi1}
		\left|(\Psi_{t}^{\lambda})''(\xi)\right|&\le C\left(\frac 1{\epsilon}+\lambda\right)e^{c\left(\frac 1{\epsilon}+\lambda\right)t}\left(1+|\xi|^{4m-1}\right),\\\label{prop-psi2}
		& \left|\Psi_{t}^{\lambda}(\xi)-\Psi_{0}^{\lambda}(\xi)\right| \le  C\left(\frac 1{\epsilon}+\lambda\right)^2e^{c\left(\frac 1{\epsilon}+\lambda\right)t} t\left(1+|\xi|^{4m+1}\right),
	\end{align}
	where $c,C>0$, $\Psi_0^{\lambda}(\xi):=-\frac 1{\epsilon} f(\xi)+\lambda \xi.$
\end{lm}

For the proof of Lemma \ref{flow-local-est}, we refer to \cite[Lemma 4.2]{CH18}. 
\revred{We also need the following a priori estimates for the scalar ODE flow in order to analyze the regularity of \eqref{spl-sch}. The proof of \eqref{comparasion} is similar to that of \cite[Lemma 1.2.6]{Cer01} and is omitted.}
\revblue{The proof of \eqref{difference} is given in Appendix \ref{appendix-difference} for completeness.} 

\begin{lm}[Dissipative flow estimates]\label{flow-dissipative-est}
Let $\lambda\ge 0$. There exists $b_0\in (0,1)$ and $C>0$ such that for every $t>0$ and $\xi\in \mathbb R$,
	\begin{align}\label{comparasion}
		{|\Phi_t^{\lambda}(\xi)|\le  C\min\left(\left(1+\left(\frac \epsilon {t}\right)^{{\frac 1{2m}}}\right), e^{-\frac {b_0}{\epsilon} t} |\xi|+ 1\right),}
	\end{align}
	Furthermore, there exists $b_1\in (0,1)$ and $C'>0$ such that for any $\xi, \xi_0\in \mathbb R$ and $t\ge 0$,
	\begin{align}\label{difference}
	{|\Phi_t^\lambda(\xi+\xi_0)-\xi_0|}
&\le
 e^{-\frac {b_1}\epsilon t}|\xi|  
+
C
\left(1+|\xi_0|^{2m+1}\right)(1-e^{-\frac {b_1}{\epsilon} t}).
\end{align}
\end{lm}

%The precise assumptions on \(A\), \(f\), and \(Q\) are collected in
%Section \ref{sec-2}.

\revred{These flow estimates will be used repeatedly in the moment, regularity, and weak error estimates. We can now state the second main result, namely the weak convergence rate of \eqref{spl-sch} with polynomial dependence on \(\epsilon^{-1}\).}

%In practice, we may take $\lambda=0$ for Dirichlet boundary condition and $\lambda>0$ such that the smallest eigenvalue of $A-\lambda$ is strictly negative.
\iffalse One can also extend our approach to analyze other temporal splitting schemes such as
\begin{align*}
	u_{n+1}=T_{\lambda,\tau}\Phi_{\tau}^{\lambda}(u_n)+T_{\lambda,\tau}\delta W_n(t),
\end{align*}
where $T_{\tau,{\lambda}}=(I-(A-\lambda)\tau)^{-1}.$ \fi

\begin{tm}\label{main-tm}
	Let \(\{u_k\}_{k\in\mathbb N}\) be the solution defined by \eqref{spl-sch}. Assume that \(x\in\mathbb H^{\beta}\cap E\), \(\beta\in(0,2)\),
	\[
	\left\|(-A+\lambda)^{\frac{\beta-1}{2}}Q^{\frac12}\right\|_{HS(\mathbb H)}<\infty,
	\]
	and that the non-degeneracy quantity \(\gamma(\epsilon)\) in
	\eqref{eq:Neps} is finite. Let \(T=N\tau\), and let
	\(\phi\in \mathcal B_{V,q}(\mathbb H)\cap\mathcal C_b^1(\mathbb H)\),
	\(q\in\mathbb N\). Then   for any \(\alpha_3\in(0,1)\), there exist a positive polynomial
	\(b_1(\|u(0)\|_E,\|u(0)\|_{\mathbb H^{\beta}}, \epsilon^{-1},\mathfrak q_\beta,\gamma(\epsilon))\) and $c_1\ge 0$ such that
	\begin{align}\label{main-tm1}
		&\left|\E\phi(u(t_N))-\E\phi(u_N)\right|\\\nonumber
		&\le
		(1+T)(\|\phi\|_{V,q}+c_1|D\phi|_\infty)b_1(\|u(0)\|_E,\|u(0)\|_{\mathbb H^{\beta}}, \epsilon^{-1},\mathfrak q_\beta,\gamma(\epsilon)) 
		\tau^{\min\left(\beta,\,\alpha_3\right)} .
	\end{align}

 In particular, if \eqref{erg-noi} holds, then for any  
	\(\phi\in\mathcal B_{V,q}(\mathbb H), q\in \mathbb N\), and \(\alpha_3\in(0,1)\), 
	\begin{align}\label{main-tm2}
		&\left|\E\phi(u(t_N))-\E\phi(u_N)\right|\\\nonumber
		&\le
		(1+T)\|\phi\|_{V,q} b_1(\|u(0)\|_E,\|u(0)\|_{\mathbb H^{\beta}}, \epsilon^{-1}, \mathfrak q_\beta,\mathfrak q_{-})\tau^{\min\left(\beta,\,\alpha_3\right)}  .
\end{align}

\end{tm}

\revblue{To the best of our knowledge, Theorem \ref{main-tm} is the first weak approximation result for a temporal semi-discretization of \eqref{sac} in the near sharp-interface regime whose error constants grow at most polynomially in \(\epsilon^{-1}\), without assuming that the noise intensity vanishes as \(\epsilon\to0\). In particular, the theorem suggests that, for a prescribed accuracy, the time step \(\tau\) needs to be only polynomially small in \(\epsilon\), rather than exponentially small, when \(\epsilon\) is small. 
When \(\epsilon=1\), the same argument also gives a time-independent weak error bound for \eqref{spl-sch} under \eqref{erg-noi}, by combining the techniques in \cite{Bre14,CHS21} with exponential ergodicity.}

\revblue{The proof is written for the temporal approximation \eqref{spl-sch}, but the analytic ingredients are separated from the specific discretization.  For instance, the same strategy can be applied to the splitting scheme}
\begin{align}\label{spl-sch1}
	u_{n+1}=S_{\lambda}(\tau)\Phi_{\tau}^{\lambda}(u_n)+S_{\lambda}(\tau)\delta W_n, \; \delta W_n=W(t_{n+1})-W(t_n),
\end{align}
\revblue{whose structure is close to \eqref{spl-sch}.  The difference is that \eqref{spl-sch} integrates the stochastic convolution exactly, while \eqref{spl-sch1} replaces it by \(S_\lambda(\tau)\delta W_n\).  Consequently, the proof for \eqref{spl-sch1} involves the second-order Kolmogorov estimate \eqref{time-ind2} to control the additional stochastic-convolution error, while the remaining terms are handled as in the proof of Theorem \ref{main-tm}; see also \cite{JKW11,CHS21}.}

\revblue{The rest of the paper is organized as follows.  Section \ref{sec-2} collects regularity estimates for the exact solution and the splitting approximation.  Section \ref{sec-4} proves the Kolmogorov regularity estimate in Theorem \ref{main-kol-reg}.  Section \ref{sec-5} proves Theorem \ref{main-tm}.  Appendix \ref{app-auxiliary-proofs} contains some auxiliary estimates.}

\section{Regularity estimates of the exact and numerical solutions}
\label{sec-2}

\revblue{This section collects the estimates for the exact dynamics and the numerical dynamics that are used later in the Kolmogorov and weak error arguments.  The estimates are stated with explicit dependence on \(\epsilon^{-1}\) and on the covariance quantity \(\mathfrak q_\beta\).}

\subsection{Properties of the original problem}

	\revblue{The existence and uniqueness of the mild solution of \eqref{sac} follow from standard arguments for monotone polynomial reaction terms; see, for example, \cite[Chapter 6]{Cer01}.  Below we collect the a priori and regularity estimates used in the rest of the paper.  The initial-value-free \(E\)-estimate in Lemma \ref{exa}(ii), which is needed for the time-independent Kolmogorov bounds, is also proved in Appendix \ref{app-initial-free-E}.}

%An alternative approach to prove Lemma \ref{exa} is to use the regularized problem in Appendix \ref{app-regularized} and then take the limit \(\delta\to0\). \rev{The initial-value-free part of Lemma \ref{exa}(i) is proved in Appendix \ref{app-initial-free-E}.}

\begin{lm}\label{exa}
	Let $u(0)\in \mathbb H^{\beta}\cap E,$ and $\left\|(-A+\lambda)^{\frac {\beta-1}2}Q^{\frac 12}\right\|_{HS(\mathbb H)}<\infty$ with $\beta\in [0,2)$ and $\lambda>0$. 
	\rev{For every \(p\ge1\), the following estimates hold:}
\begin{enumerate}
\item[(i)]  \rev{There exist \(c_1\in(0,1)\) and a positive polynomial
\(P_{2m+1}(\mathfrak q_\beta)\) of degree $2m+1$ such that for every $t\ge0,$}
\begin{align}\label{eq:initial-free-E0}
 \|u(t)\|_{L^p(\Omega;E)}
  \le e^{-c_1t/\epsilon} \|u(0)\|_{E}+\rev{P_{2m+1}(\mathfrak q_\beta)}.
\end{align}
\rev{Moreover, for any $p\ge 1$, there exists a positive polynomial
\(P_{p}(\mathfrak q_\beta)\)  of order $p$  such that for any $t>0$,}
\begin{align}\label{eq:initial-free-E}
\|u(t)\|_{L^p(\Omega;E)}
\le
P_{p}(\mathfrak q_\beta)\left(1+\left(\frac {\epsilon}t\right)^{\frac 1{2m}} \right).
\end{align}
\item[(ii)]For any \(p\ge 1\), there exists a positive polynomial
\(P_{(2m+1)^2,2m+1,1}(\mathfrak q_\beta, \|u_0\|_E$, $\|u_0\|_{\mathbb H^{\beta}})\)
\revred{of degrees \((2m+1)^2\), \(2m+1\), and \(1\) in its three variables, respectively, such that for any \(t\ge0\),}
\begin{align}\label{eq:initial-free-E1}
  \|u(t)\|_{L^p(\Omega;\HH^\beta)}
  &\le \rev{P_{(2m+1)^2,2m+1,1}(\mathfrak q_\beta, \|u_0\|_E,\|u_0\|_{\mathbb H^{\beta}})}\frac 1\epsilon.
\end{align}
\revred{and for every \(s,t\ge0\),}
\begin{align}\label{eq:initial-free-E2}
 \|u(t)-u(s)\|_{L^p(\Omega;\HH)}
  &\le \rev{ P_{(2m+1)^2,2m+1,1}(\mathfrak q_\beta,\|u_0\|_E,\|u_0\|_{\mathbb H^{\beta}})}  
\frac 1\epsilon
  |t-s|^{\min\{\beta/2,1/2\}} .
\end{align}

\end{enumerate}

	\iffalse
	  It holds that for some $c_0\in (0,1)$ and $C(p,Q)>0$,
	\begin{align*}
		\|u(t)\|_{L^{p}(\Omega; E)} \le e^{-\frac {c_0 t}{\epsilon} } \|u(0)\|_{L^{p}(\Omega; E)}+C(p,Q),
	\end{align*}
	Moreover, when $\beta<2$, we have that for $\gamma\in [0,\beta]$ and $C({p,Q,\beta,\gamma})>0,$
	\begin{align*}
		\left\|u(t)\right\|_{L^{p}(\Omega; \mathbb H^{\gamma})} &\le C({p,Q,\beta,\gamma})\left(\lambda+\frac 1\epsilon\right)^{\frac {\gamma} 2}\left(1+\|u(0)\|_{\mathbb H^{\beta}}+\|u(0)\|_{E}^{2m+1}\right)^{\frac {\gamma}{\beta}}\\
		&\quad \times\left(1+\|u(0)\|\right)^{1-\frac {\gamma}{\beta}}
	\end{align*}
	and that for $\theta\in [0,1]$ and $C(p,Q,u(0),\beta,\theta)>0,$ $s,t\ge0,$
	\begin{align*}
		\|u(t)-u(s)\|_{L^{p}(\Omega; \mathbb H)} \le C(p,Q,u(0),\beta,\theta) \left(\lambda+\frac 1 {\epsilon}\right)^{\beta \theta} \left|t-s\right|^{\min\left(\frac {\beta \theta} 2,\frac 12\right)}.
	\end{align*}
	When $\beta=2$, it holds that for some $C({p,Q,u(0)})>0,$ $s,t\ge0,$
	\begin{align*}
		\|u(t)\|_{L^{p}(\Omega; \mathbb H^{2})} 
		&\le C({p,Q,u(0)})\left(\frac 1{\epsilon}+\lambda\right),\\
		\|u(t)-u(s)\|_{L^p(\Omega;\mathbb H)}&\le C({p,Q,u(0)})\left(\frac 1{\epsilon}+\lambda\right)\left|t-s\right|^{\frac 12}.
	\end{align*}
	\fi
\end{lm}

\subsection{A priori bound  of splitting numerical approximation}

In this part, we establish an a priori bound for the splitting approximation that is independent of $\epsilon$. As a by-product, the regularity estimates of the splitting scheme are shown to depend at most polynomially on $ \frac {1} {\epsilon} $.  Note that our approach to tracking the polynomial dependence on \(\epsilon^{-1}\) in the weak error analysis can also be extended to other temporal and spatial discretizations, provided that they satisfy similar time-independent moment bounds in \(E\) and \(\mathbb H^\beta\).

%In this section, we introduce a concrete numerical approximation considered in the rest of the paper. It can be seen that our approach to studying the polynomial dependence on $\frac 1\epsilon$  in weak error analysis could be extended to deal with other temporal and spatial discretizations as long as they enjoy similar time-independent moment bounds in $E$ and $\mathbb H^{\gamma}, \gamma \in (0,\beta)$. 

\begin{lm}\label{pri-un-e}
	Let $u(0)\in \mathbb H^{\beta}\cap E,$ and $\Big\|(-A+\lambda)^{\frac {\beta-1}2}Q^{\frac 12}\Big\|_{HS(\mathbb H)}<\infty$ with $\beta\le 2$. 
	For any $p\ge 1$, there exists a positive polynomial $P_{2m+1}(\mathfrak q_\beta)$ of order $2m+1$ such that 
	\begin{align*}
		\rev{\sup_{n\in \mathbb N} \left \|u_{n} \right\|_{L^p(\Omega;E)}
		\le P_{2m+1}(\mathfrak q_\beta)(1+\|u_0\|_E)}.
	\end{align*}
\end{lm}

\begin{proof}
Denote the stochastic convolution at the time grid $t_{n+1}$ by 
\[
 z_{n+1}=S_{\lambda}({\tau}) z_n+S_{\lambda}({\tau})\delta W_n, \; z_0=0.
\]
Decompose $u_{n+1}=S_{\lambda}({\tau})\Phi_\tau^\lambda(u_n)+S_{\lambda}({\tau})\delta W_n$ by \[
 u_{n+1}=v_{n+1}+z_{n+1}.
\]
Then 
we have
\begin{equation}\label{eq:v-recursion}
 v_{n+1}
 =S_{\lambda}({\tau})\bigl(\Phi_\tau^\lambda(v_n+z_n)-z_n\bigr).
\end{equation}

	First, since $z_{n+1}=z(t_{n+1})$, by \cite[Lemma 8.2.1 and Proposition 8.2.2]{Cer01}, it holds that 
	\begin{align}\label{pri-zn}
		&\sup_{n\in \mathbb N}\left\|z_{n}\right\|_{L^p(\Omega;E)}+\sup_{n\in \mathbb N}\left\|z_{n}\right\|_{L^p(\Omega; \mathbb H^{\beta})}\le C_{p,\beta} \mathfrak q_\beta.
	\end{align}

	\iffalse
	Denote $v_n=u_n-z_n$. Then by choosing suitable $\lambda>0$, applying the contraction property $\|S_{\lambda}(\tau) v\|_{E}\le e^{-\lambda t} \|v\|_E$, Young's and H\"older's inequalities, it follows that 
	\begin{align*}
		\left\|v_{n+1}\right\|_E&\le \left\|S_{\lambda}({\tau})(\Phi_{\tau}^{\lambda} \left(v_n+z_n)-z_n\right)\right\|_E
		\\
		&\le e^{-\lambda \tau}\left\|\Phi_{\tau}^{\lambda} (v_n+z_n)-z_n\right \|_E.
	\end{align*}
%	$ \le e^{-\lambda\tau-b_0\tau}\left\|v_n\right\|_E+e^{-\lambda\tau-b_0\tau} \frac {b_{1,n}}{b_0}\left(1-e^{-b_0\tau}\right),$
\fi
	
 By the contraction property of the semigroup $S_\lambda(\cdot)$,
\[
\|S_{\lambda}({\tau})\|_{\mathcal L(E)}\le 1,
\]
we infer from \eqref{difference} that for some $b_0\in (0,1),$
\begin{align}\label{eq:pathwise-v}
\|v_{n+1}\|_E
&\le
 e^{-\frac {b_0} \epsilon \tau}\|v_n\|_E+
C
\left(1+\|z_n\|_E^{2m+1}\right)(1-e^{-\frac {b_0} \epsilon \tau}).
\end{align}
Thus,
\begin{equation}\label{eq:v-discrete}
\|v_{n+1}\|_E
\le
\rho\|v_n\|_E
+
C(1-\rho)
\left(1+\|z_n\|_E^{2m+1}\right),
\qquad
\rho=e^{-b_0\tau/\epsilon}\in(0,1).
\end{equation}
Iterating \eqref{eq:v-discrete} yields
\begin{equation}\label{eq:v-iterate}
\|v_{n+1}\|_E
\le
\rho^{n+1}\|u_0\|_E
+
C\sum_{k=0}^n \rho^{n-k}(1-\rho)
\left(1+\|z_k\|_E^{2m+1}\right).
\end{equation}
Taking the $L^p(\Omega)$-norm and using Minkowski's inequality gives
\begin{align*}
\|v_{n+1}\|_{L^p(\Omega;E)}
&\le
\|u_0\|_E
+
C\sum_{k=0}^n \rho^{n-k}(1-\rho)
\left(1+
\|\|z_k\|_E^{2m+1}\|_{L^p(\Omega)}
\right)
\\
&\le
\|u_0\|_E
+
C\left(1+
\sup_{k\in\mathbb N}\|z_k\|_{L^{p(2m+1)}(\Omega;E)}^{2m+1}
\right)
\sum_{k=0}^n \rho^{n-k}(1-\rho).
\end{align*}
Since
\[
\sum_{k=0}^n \rho^{n-k}(1-\rho)
\le1,
\]
we obtain, by \eqref{pri-zn}, there exists a positive polynomial $P_{2m+1}(\mathfrak q_\beta)$ of degree $2m+1$ such that 
\begin{equation}\label{eq:v-final}
\sup_{n\in\mathbb N}\|v_n\|_{L^p(\Omega;E)}
\le
P_{2m+1}(\mathfrak q_\beta)\rev{(1+\|u_0\|_E).}
\end{equation}
Finally,
\[
\|u_n\|_E\le \|v_n\|_E+\|z_n\|_E.
\]
Combining \eqref{eq:v-final} with the stochastic convolution bound
\eqref{pri-zn} completes the proof.

\end{proof}

\begin{rk}\label{rk:moment-other-splitting}
The same \(E\)-moment estimate holds for the alternative splitting scheme
\eqref{spl-sch1}.  In that case only the estimate of discrete stochastic convolution is
changed, and the proof above applies with the corresponding Ornstein--Uhlenbeck
recursion for \(z_n\).
\end{rk}

To end this section, we present the temporal and spatial regularity estimates of the discrete scheme. 
\begin{prop}\label{dis-con}
	Let $u(0)\in \mathbb H^{\beta}\cap E,$ and $\Big\|(-A+\lambda)^{\frac {\beta-1}2}Q^{\frac 12}\Big\|_{HS(\mathbb H)}<\infty$ with $\beta\in (0,2)$ and $\lambda>0$.  Let $\tau\lesssim \epsilon.$
	For any $p\ge 1$,  there exists a positive polynomial $P_{(2m+1)^2,2m+1,1}(\mathfrak q_\beta,\|u(0)\|_E,\|u(0)\|_{\mathbb H^{\beta}})$ of orders $(2m+1)^2,2m+1$ and 1, with respect to variables $\mathfrak q_\beta,\|u(0)\|_E,\|u(0)\|_{\mathbb H^{\beta}}$, such that any $n,k\le N,$
	\begin{align}\label{3.8}
		\|u_{k}-u_{n}\|_{L^p(\Omega;\mathbb H)}&\le \rev{P_{(2m+1)^2,2m+1}(\mathfrak q_\beta,\|u(0)\|_E,\|u(0)\|_{\mathbb H^{\beta}})}\frac {1}\epsilon |t_n-t_k|^{\min \rev{\left\{\frac {\beta} 2,\frac1 2\right\}}}.
	\end{align}
	\iffalse When $\beta=2$, 
	letting $\tau \sim O\left(\left(\frac 1\epsilon+\lambda\right)^{-1}\right)$, it holds that for $\theta\in (0,1],$
	\begin{align}\label{3.9}
		\|u_{n+1}-u_n\|_{L^p(\Omega;\mathbb H)}&\le C({p,Q,u(0)},\theta)\left(\frac1\epsilon+ \lambda\right)^{1+\theta}(1+\lambda)^{-\theta }\tau^{\min(\theta,\frac 12)}.
	\end{align}\fi
	and that 
		\begin{align}\label{reg-nu}
			\sup_{\rev{n\le N}}\left\|u_n\right\|_{L^{p}(\Omega; \mathbb H^{\beta})}& \le \rev{P_{(2m+1)^2,2m+1,1}(\mathfrak q_\beta,\|u(0)\|_E,\|u(0)\|_{\mathbb H^{\beta}})}\frac 1\epsilon.
		\end{align}
\end{prop}

\begin{proof}
	According to the regularity estimates of $z_n$ \eqref{pri-zn} and \eqref{hold-convolution}, it suffices to derive those of $v_n$.
	%Thanks to the equivalence of the norms between the fractional Sobolev space  $W^{\gamma,2}$ and the interpolation space $\mathbb H^{\gamma}$. In particular,   the function in $\mathbb H^{\gamma}$ is considered as an element in $W^{\gamma, 2}$ under the Dirichlet boundary condition with $\gamma>\frac 12$ or under the Neumann boundary condition with $\gamma>\frac 32$. 
	From the analytical properties of $S_{\lambda}({\tau})$ \eqref{smooth0}, it follows that for  some $c_0>0,$ 
	\begin{align*}
		&\|v_{n+1}\|_{\mathbb H^{\beta}}\\
		&\le \left\|(S_{\lambda}({\tau}))^{n+1}v_0\right\|_{\mathbb H^{\beta}}+\Big\|\sum_{k=0}^n(S_{\lambda}({\tau}))^{n+1-k} \Psi^{\lambda}_{\tau}(v_k+z_k) \tau\Big\|_{\mathbb H^{\beta}}\\
		&\le e^{-c_0t_{n+1}} \|v_0\|_{\mathbb H^{\beta}}
		+C\sum_{k=0}^n e^{-c_0 (t_{n+1}-t_k)} \left(t_{n+1}-t_k\right)^{-\frac {\beta}2} \left\|\Psi_{\tau}^{\lambda}(v_k+z_k) \right\|\tau.
	\end{align*}
Since $\tau\lesssim \epsilon$, applying Lemma \ref{pri-un-e}, \eqref{prop-psi} and \eqref{pri-zn},
	we have that 	\begin{align*}
		&\|v_{n+1}\|_{L^{p}(\Omega;\mathbb H^{\beta})}\\
		&\le e^{-c_0 t_{n+1}} \|v_0\|_{\mathbb H^{\beta}}+C\sum_{k=0}^n e^{-c_0\lambda (t_{n+1}-t_k)} (t_{n+1}-t_k)^{-\frac {\beta}2} \tau\left(\frac  1{\epsilon}+\lambda\right)\\
		& \quad \times \sup_{k\le n} \left(1+\|v_k\|_{{L^{p}(\Omega; L^{4m+2})}}^{2m+1}+\|z_k\|_{{L^{p}(\Omega; L^{4m+2})}}^{2m+1}\right)\\
		&\le e^{-c_0 t_{n+1}} \|v_0\|_{\mathbb H^{\beta}}+C\frac1\epsilon\sup_{k\le n}\left(1+\|v_k\|_{L^p(\Omega;E)}^{2m+1}+\|z_k\|_{L^p(\Omega;E)}^{2m+1}\right)\\
		&\le C\frac 1\epsilon\left(1+\|u_0\|_{\mathbb H^{\beta}}+\|u_0\|_{E}^{2m+1}+\mathfrak q_{\beta}^{(2m+1)^2}\right),
	\end{align*}
	where in the second inequality, we also use the fact that 
	\begin{align*}
    \tau\sum_{k=0}^{n}e^{-c_0(t_{n+1}-t_k)}
    (t_{n+1}-t_k)^{-\beta/2}
    &\le C\int_0^\infty e^{-c_0 r}r^{-\beta/2}\,dr<+\infty.
\end{align*}
The above estimate implies \eqref{reg-nu} by taking the supremum over $n$.

Next, we prove the discrete H\"older continuity estimate \eqref{3.8}. For simplicity, we assume that $n>k.$
Let \(h=t_n-t_k\).  From the mild representation of the scheme,
\[
    v_n=S_\lambda(t_k)v_k+\tau\sum_{j=k}^{n-1}
    S_\lambda(t_n-t_j)\Psi_\tau^\lambda(u_j),
\]
and hence
\begin{align*}
    v_n-v_k
    &= (S_\lambda(h)-I)v_k
    +\tau\sum_{j=k}^{n-1}S_\lambda(t_n-t_j)\Psi_\tau^\lambda(u_j).
\end{align*}
We estimate the above two terms separately.

First, for \(\beta\in (0,2)\), the analytic semigroup estimate \eqref{smooth1} with $\kappa_1=\frac \beta 2$,  and \eqref{reg-nu} give
\begin{align*}
    \|{(S_\lambda(h)-I)v_k}\|_{L^p(\Omega;\mathbb H)}
    &\le C h^{\beta/2}\|v_k\|_{L^p(\Omega;\mathbb H^{\beta})}\\
    &\le Ch^{\beta/2}\frac 1\epsilon P_{(2m+1)^2,2m+1,1}(\mathfrak q_\beta,\|u(0)\|_E,\|u(0)\|_{\mathbb H^{\beta}}).
\end{align*}
Second, by  \eqref{smooth0} with $\kappa=0$, using Lemma \ref{pri-un-e}, \eqref{prop-psi}, and \eqref{pri-zn}, we have 
\begin{align*}
    &\tau\sum_{j=k}^{n-1}\|{S_\lambda(t_n-t_j)\Psi_\tau^\lambda(u_j)}\|_{L^p(\Omega;\mathbb H)}\\
    &\le C \frac 1\epsilon \sum_{j=k}^{n-1} e^{-c_0(t_n-t_j)} \tau \sup_{k\le n} \left(1+\|v_k\|_{{L^{p}(\Omega; L^{4m+2})}}^{2m+1}+\|z_k\|_{{L^{p}(\Omega; L^{4m+2})}}^{2m+1}\right)
   \\
   &\le  C \frac 1\epsilon \min(h,1) \left(1+\|u(0)\|_{E}^{2m+1}+\mathfrak q_\beta^{(2m+1)^2}\right).
\end{align*}
Combining the two estimates above with \eqref{hold-convolution}, we conclude that \eqref{3.8} holds.

\end{proof}

\iffalse
\begin{proof}
	By decomposing  $\|u_{n+1}-u_n\|_{L^p(\Omega;\mathbb H)}$ into $\|v_{n+1}-v_n\|_{L^p(\Omega;\mathbb H)}$ and $\|z_{n+1}-z_n\|_{L^p(\Omega;\mathbb H)}$, it holds that 
	\begin{align*}
		\|z_{n+1}-z_n\|_{L^p(\Omega;\mathbb H)}\le \|(S_{\tau}^{\delta}-I)z_n\|_{L^p(\Omega;\mathbb H)}+\|S_{\tau}^{\delta}\delta_n W\|_{L^p(\Omega;\mathbb H)}\le C \tau^{\beta}.
	\end{align*}
	It remains to estimate $\|v_{n+1}-v_n\|_{L^p(\Omega;\mathbb H)}$. According to Proposition \ref{reg-nu} and the a priori estimate of $z_n$, we obtain
	\begin{align*}
		\|v_{n+1}-v_n\|_{L^p(\Omega;\mathbb H)}&\le \|(S_{\lambda}({\tau})-I)v_n\|_{L^p(\Omega;\mathbb H)}+\tau \|S_{\lambda}({\tau}) \Psi_{\tau}^{\lambda}(v_n+z_n)\|_{L^p(\Omega;\mathbb H)}\\
		&\le C(\lambda,m,Q,p,u_0,\gamma)[(1+\frac 1 {\epsilon})^{\frac {\gamma}{\beta}} \tau^{\gamma}+\frac {\tau}{\epsilon}].\end{align*}\end{proof}
\fi

\section{Regularity estimates for the Kolmogorov equation}
\label{sec-4}
\revblue{In this section, we prove the time-independent regularity estimate of the
Kolmogorov equation stated in Theorem \ref{main-kol-reg}. }

\subsection{Regularity estimates of the Kolmogorov equation in a bounded interval}

\revred{We first prove finite-time regularity estimates for the original Kolmogorov equation. These estimates may depend exponentially on \(\epsilon^{-1}\) and on the time horizon \(T\).}

\begin{prop}\label{reg-kl-short}
	Let $x\in \mathbb H^{\beta}\cap E,$ and $\Big\|\left(-A+\lambda\right)^{\frac {\beta-1}2}Q^{\frac 12}\Big\|_{HS(\mathbb H)}<\infty$ with $\beta<2$.
	Let $T>0$ and $\phi\in \mathcal C_b^2(\mathbb H).$
	For $\alpha\in (\frac d4,1),\alpha_1,\alpha_2\in [0,1)$, $\alpha_1+\alpha_2<1,$ there exist a positive constant $C'$ and positive polynomials \(P_{2m(2m+1),2m}(\mathfrak q_\beta,\|x\|_E)\), \(P_{(8m-1)(2m+1),8m-1}(\mathfrak q_\beta,\|x\|_E)\),  such that for any $h,h_1,h_2\in H$ and $t\in (0,T]$,
	\begin{align}\label{reg-kol-1}
		|DX(t,x)\cdot h|&\le |D\phi|_\infty P_{2m(2m+1),2m}(\mathfrak q_\beta,\|x\|_E)\frac 1{\epsilon}\exp\Big(\frac {C't}{\epsilon} \Big)\\\nonumber 
		&\quad  \times \revred{\| h\|_{\mathbb H^{-2\alpha_1}}}
		(1+t^{-\alpha_1}),\\\nonumber
		|D^2 X(t,x)\cdot (h_1,h_2)|&\le (|D\phi|_\infty +|D^2\phi|_\infty )P_{(8m-1)(2m+1),8m-1}(\mathfrak q_\beta,\|x\|_E)  \left(\frac1\epsilon\right)^{4} 
		\\ \nonumber 
		&\quad\times \exp\Big(\frac {C't}{\epsilon} \Big)\Big(t^{1-\alpha-\alpha_1-\alpha_2}+t^{-\alpha_1-\alpha_2}
		+1\Big) \\ \label{reg-kol-2}
		&\quad\times \revred{\|h_1\|_{\mathbb H^{-2\alpha_1}} 
		\|h_2\|_{\mathbb H^{-2\alpha_2}}}.
	\end{align}
\end{prop}

 \revblue{Since the proof is close to \cite[Proposition 4.1]{CH18}, we only use the result here and give the details in Appendix \ref{appendix-reg-finite} for completeness.}
Later, we will improve this exponential dependence using the asymptotically strong Feller estimate or the strong Feller estimate.

\subsection{Longtime regularity estimate of Kolmogorov equation}

\revblue{Our strategy for overcoming the exponential dependence on \(\epsilon^{-1}\) is to combine the short-time estimates above with the long-time smoothing of the Markov semigroup \(P_t\).  
We use the weighted spaces \(\mathcal B_{V,q}(\mathbb H)\) introduced in
Section \ref{sec-main}.}

Recall the low-mode non-degeneracy quantity \(\gamma(\epsilon)\) defined in
\eqref{eq:Neps}.  In this section we assume \(\gamma(\epsilon)<\infty\).  The
non-degeneracy condition \eqref{erg-noi} is used only when a terminal coupling or a
second-order derivative estimate is needed.

\iffalse
we assume that 
\begin{align}\label{erg-noi}
	\rev{\left\|Q^{-\frac 12}(-A+\lambda I)^{-\frac 12}\right\|_{\mathcal L(\mathbb H)}<\infty.}
\end{align}
\rev{This assumption is the non-degeneracy condition used in \cite{CHS21}; the norm is the operator norm on \(\mathbb H\).  If \(Qe_j=q_je_j\) and \(-Ae_j=\lambda_je_j\), then \eqref{erg-noi} is equivalent to}
\[
\rev{\sup_{j\ge1}q_j^{-1/2}(\lambda_j+\lambda)^{-1/2}<\infty.}
\]
\rev{It is compatible with the Hilbert--Schmidt condition in Theorem \ref{main-tm} for a range of trace-class noises.  For example, if \(q_j\asymp(\lambda_j+\lambda)^{-\rho}\), then \eqref{erg-noi} requires \(\rho\le1\), while \(\|(-A+\lambda)^{(\beta-1)/2}Q^{1/2}\|_{HS(\mathbb H)}<\infty\) requires \(\rho>\beta-1+d/2\).  Hence, admissible examples exist whenever \(\beta<2-d/2\).  This compatibility restriction is one reason why the total variation statement below is separated from the smoother weak-error estimates.}
\fi

\revblue{It is known that the full non-degeneracy condition \eqref{erg-noi} implies strong Feller and irreducibility properties for the Markov semigroup; see \cite{CHS21}.  The estimate below is a quantitative version adapted to the sharp-interface regime.  It is stated first under the low-mode condition \(\gamma(\epsilon)<\infty\) and then under the stronger condition \eqref{erg-noi}.  The explicit polynomial dependence will be used in the proof of Theorem \ref{main-tm}.}

\begin{lm}\label{Strong Feller property}
	Let $x\in \mathbb H^{\beta}\cap E,$ and $\Big\|\left(-A+\lambda\right)^{\frac {\beta-1}2}Q^{\frac 12}\Big\|_{HS(\mathbb H)}<\infty$ with $\beta<2$. 	Let $\phi\in \mathcal B_{V,q}(\mathbb H)$ with $q\in \mathbb N.$ Let $\gamma(\epsilon)<+\infty.$
		\revblue{Then there exist \(c_0>0\), \(c_1\ge0\), and a positive function \(\tilde b(\epsilon^{-1},t,\mathfrak q_\beta,\gamma(\epsilon))\), polynomial in its displayed non-time variables, such that for any \(t>0\),}
	\begin{align}\label{asy-strong-feller}
	\|DP_t\phi\|_{\infty}\le \tilde b(\epsilon^{-1},t, \mathfrak q_\beta,\gamma(\epsilon))\|\phi\|_{V,q}+c_1e^{-\frac {c_0}\epsilon t}\|D\phi\|_{\infty}.	
	\end{align}
	Here $\tilde b(\epsilon^{-1},t,\mathfrak q_\beta,\gamma(\epsilon))=\frac{(1+L_f)\gamma(\epsilon)}{\sqrt{\epsilon}}
P_{q^2}(\mathfrak q_{\beta})
\left(1+\left(\frac{\epsilon}{t}\right)^{q/(2m)}\right)$ for some polynomial $P_{q^2}(\cdot)$ of degree $q^2$.
	
		\revblue{In particular, if \eqref{erg-noi} holds, then \eqref{asy-strong-feller} holds with \(c_1=0\) and with \(\tilde b(\epsilon^{-1},t, \mathfrak q_\beta$, $\gamma(\epsilon))\) replaced by \(\tilde b(\epsilon^{-1},t,\mathfrak q_\beta,\mathfrak q_-)\).} 	
\end{lm}

\begin{proof}

For \eqref{asy-strong-feller} under the condition $\gamma(\epsilon)<+\infty$, we provide the proof in Appendix \ref{general-time-ind}, since it follows the asymptotic coupling argument in the proof of \cite[Lemma 4.1]{CW2024}. Here we present the detailed proof under the nondegeneracy condition \eqref{erg-noi}.

We first give the standard coupling-by-change-of-measure proof using a terminal
coupling control. \revblue{The argument can be justified by applying it to Galerkin approximations of \eqref{sac}, deriving estimates independent of the Galerkin dimension, and then passing to the limit.  We use this standard approximation procedure implicitly below.}

Fix $t>0.$ 
Let $u(t,x)$ be the solution of \eqref{sac} starting from $x$.
Define a controlled process $Y$ with $Y(0)=y$ by
\begin{equation}\label{controlled-Y}
    dY(s)=AY(s)\,ds\rev{-\epsilon^{-1}f(Y(s))\,ds}+Q^{\frac 12}v(s)\,ds+dW(s),
    \qquad 0\le s\le t.
\end{equation}
Set
\[
    \tilde Z(s)=u(s)-Y(s),
    \qquad
    \tau_{sto}=\inf\{s\ge0: \tilde Z(s)=0\}.
\]
For $s<\tau_{sto}$, choose the control through
\begin{equation}\label{correct-control}
    Q^{\frac 12}v(s)
    =
    \gamma(s)\frac{\tilde Z(s)}{\|\tilde Z(s)\|},
    \qquad
    \gamma(s)
    =
    \frac{\|{x-y}\|e^{-\frac {L_f}{\epsilon} s}}
    {\int_0^t e^{-2\frac {L_f}{\epsilon} r}\,dr},
\end{equation}
and set $v(s)=0$ for \(\rev{s\ge\tau_{sto}}\). 
Here, $L_f>0$ is the one-sided Lipschitz constant of $f$. Notice that the control is proportional to
$\tilde Z(s)/\|{\tilde Z(s)}\|$. This is the point that makes the
coupling occur exactly at the terminal time.

For $s<\tau_{sto}$, using one-sided Lipschitz continuity of $f$ and integration by parts, we obtain
\begin{align}\label{coupling-est}
    \frac12\frac{d}{ds}\|{\tilde Z(s)}\|^2+\|\nabla \tilde Z(s)\|^2    \le
   \frac {L_f}\epsilon \|{\tilde Z(s)}\|^2-
    \gamma(s)\|{\tilde Z(s)}\|.
\end{align}
Omitting the term $\|\nabla \tilde Z(s)\|^2$, and dividing by $\|\tilde Z(s)\|$ gives
\begin{equation}\label{distance-ineq}
    \frac{d}{ds}\|\tilde Z(s)\|
    \le
   \frac {L_f}\epsilon\|\tilde Z(s)\|-\gamma(s),
    \qquad s<\tau_{sto}.
\end{equation}
Multiplying \eqref{distance-ineq} by $e^{-\frac {L_f}{\epsilon} s}$ yields
\[
    \frac{d}{ds}\left(e^{-\frac {L_f}{\epsilon} s}\|{\tilde Z(s)}\|\right)
    \le
    -e^{-\frac {L_f}{\epsilon} s}\gamma(s)
    =
    -\frac{\|{x-y}\|e^{-2\frac {L_f}{\epsilon} s}}
    {\int_0^t e^{-2\frac {L_f}{\epsilon}  r}\,dr}.
\]
Hence, for $s\le t\wedge\tau_{sto}$,
\[
    e^{-\frac {L_f}{\epsilon} s}\|{\tilde Z(s)}\|
    \le
    \|{x-y}\|
    \left(
    1-
    \frac{\int_0^s e^{-2\frac {L_f}{\epsilon} r}\,dr}
    {\int_0^t e^{-2\frac {L_f}{\epsilon}r}\,dr}
    \right).
\]
At $s=t$ the right-hand side is zero. Therefore $\tau_{sto}\le t$ and
\[
    u(t)=Y(t).
\]
This proves terminal coupling.

Thanks to \eqref{erg-noi}, for $s<\tau_{sto}$,
\begin{align*}
    \|v(s)\|^2&= \Big\|Q^{-\frac 12}\gamma(s)\frac{\tilde Z(s)}{\| \tilde Z(s)\|} \Big\|^2
    \le
   \gamma(s)^2 \frac {\|Q^{-\frac 12}\|_{\mathcal L(H)}^2\|\tilde Z(s)\|^2}{\|  \tilde Z(s)\|^2}  \le \mathfrak q_-^2  \gamma(s)^2.
\end{align*}
Therefore
\begin{equation}\label{control-cost}
    \int_0^t\|{v(s)}\|^2\,ds
    \le
   \mathfrak q_-^2\int_0^t\gamma(s)^2\,ds
    =
    \mathfrak q_-^2\frac{\|{x-y}\|^2}{\int_0^t e^{-2\frac {L_f} \epsilon  r}\,dr}
    =
    \mathfrak q_-^2\frac{2\frac {L_f}\epsilon}{1-e^{-2\frac {L_f}\epsilon t}}
    \|{x-y}\|^2.
\end{equation}

Define the Girsanov density
\[
    \mathcal E_t
    =
    \exp\left(
    -\int_0^t \<{v(s)},{dW(s)}\>
    -\frac12\int_0^t\|{v(s)}\|^2\,ds
    \right).
\]
Under the probability measure with density $\mathcal E_t$, the controlled
process $Y$ has the same law as the uncontrolled solution starting from $y$.
Since $u(t)=Y(t)$, for positive $\phi$,
\[
    P_t\log\phi(y)
    =\E[\mathcal E_t\log\phi(Y(t))]
    =\E[\mathcal E_t\log\phi(u(t))].
\]
The entropy inequality gives
\[
    \E[\mathcal E_t\log\phi(u(t))]
    \le
    \log\E[\phi(u(t))]+
    \E[\mathcal E_t\log\mathcal E_t].
\]
Moreover,
\[
    \E[\mathcal E_t\log\mathcal E_t]
    =
    \frac12\E_{\mathcal E}\int_0^t\|{v(s)}\|^2\,ds
    \le
   \mathfrak q_-^2\frac{\frac {L_f}{\epsilon}}{1-e^{-2\frac {L_f}{\epsilon} t}}
   \|{x-y}\|^2.
\]
Combining the above estimates proves 
\begin{equation}\label{finite-logh}
    P_t\log\phi(y)
    \le
    \log P_t\phi(x)
    +
   \mathfrak q_-^2  \frac{\frac {L_f}\epsilon}{1-e^{-2\frac {L_f}\epsilon t}}
    \|{x-y}\|^2,
    \qquad t>0,
\end{equation}
This, together with the standard implication from the log-Harnack inequality to
gradient estimates, implies the desired result \eqref{asy-strong-feller} for
\(\phi\in \mathcal B_b(\mathbb H)\), $q=0$, where \(c_1=0\) and $\tilde b$ is replaced by 
\begin{align}\label{b-strong-feller}
b(\frac 1\epsilon, t,\mathfrak q_-) =\mathfrak q_- \sqrt{\frac {1}{2t}}+\mathfrak q_- \sqrt{\frac {L_f}{\epsilon}}\ge \mathfrak q_- \sqrt{\frac{\frac {L_f}\epsilon}{1-e^{-2\frac {L_f}\epsilon t}}}.
\end{align}

Next, we use the Markov property $P_t\phi=P_{\frac t2} P_{\frac t2}\phi$ to show that \eqref{asy-strong-feller} also holds for $\phi \in \mathcal B_{V,q}(\mathbb H)$ for any $q\in \mathbb N^+$ with $c_1=0.$  We first recall the \(x\)-independent a priori estimate \eqref{eq:initial-free-E} which is  needed for this extension.

Now we prove \eqref{asy-strong-feller} for $\phi \in \mathcal B_{V,q}(\mathbb H)$.
By the Markov property of $P_t$, it holds that 
\begin{align*}
P_t\phi(x)=P_{\frac t2} \E[\phi(u(\frac t2,x))]=P_{\frac t2}(P_{\frac t2}\phi)(x)
\end{align*}
Therefore,
\begin{align*}
\|DP_t\phi\|_{\infty}\le b(\frac 1\epsilon,\frac t2,\mathfrak q_{-})\|P_{\frac t 2}\phi\|_{\infty}.
\end{align*}
By \eqref{eq:initial-free-E}, 
\begin{align*}
P_{\frac t 2}\phi(x)=\E[\phi(u(\frac t2,x))]&\le C\|\phi\|_{V,q} \E[1+\|u(\frac t2,x)\|^q]\\
&\le \|\phi\|_{V,q} (P_{q}(\mathfrak q_{\beta}))^q \left(1+\left(\frac{\epsilon}{t}\right)^{q/(2m)}\right).
\end{align*}
As a consequence, 
we have 
\begin{align*}
\|DP_t\phi\|_{\infty}\le   b(\frac 1\epsilon, \frac t2,\mathfrak q_-) (P_{q}(\mathfrak q_{\beta}))^q \|\phi\|_{V,q}  \left(1+\left(\frac{\epsilon}{t}\right)^{q/(2m)}\right).
\end{align*}
This implies \eqref{asy-strong-feller} with 
\begin{align}\label{asy-strong-feller1}
\tilde b(\frac 1\epsilon,t,\mathfrak q_{\beta},\mathfrak q_-)=b(\frac 1\epsilon, \frac t2,\mathfrak q_-) P_{q^2}(\mathfrak q_{\beta}) \left(1+\left(\frac{\epsilon}{t}\right)^{q/(2m)}\right).
\end{align}

\iffalse

\begin{align*}
d\|v(t)\|^2&=2 \<Av, v\>dt-\frac 2{\epsilon}\<f(v+z),v\>dt\\
&\le -\frac 2\epsilon (a_{2m+1}-\eta)\|v\|_{L^{2m+2}}^{2m+2}dt+\frac 2{\epsilon}C(\eta,m)(1+\|z\|_{L^{2m+2}}^{2m+2}) dt
\end{align*}

	On the other hand,
	notice that Young's inequality yields that for $\eta<1$ and any $\xi\in \mathbb R,$ 
	\begin{align*}
		-\frac 1 {\epsilon}f(\xi)+\lambda \xi \le-\frac 1{\epsilon}(1-\eta) a_{2m+1}\xi^{2m+1}+ \left(\frac 1\epsilon C(\eta,m) +\lambda\right) \xi.
	\end{align*}
	\rev{Considering the parameterised Bernoulli equation}
	\begin{align*}
		dx(t)=-\frac 1{\epsilon}(1-\eta)a_{2m+1} x(t)^{2m+1}+\left(\frac 1\epsilon C(\eta,m) +\lambda\right) x(t),
	\end{align*}
	\rev{and applying the one-dimensional comparison principle, we obtain the upper bound by solving this Bernoulli equation explicitly:}  
	\begin{align*}
		\left|\Phi_t^{\lambda}(x_0)\right|^{2m}&\le \frac {|x_0|^{2m}}{e^{-2m(\frac 1\epsilon C(\eta,m) +\lambda) t}+\frac { (1-\eta)a_{2m+1}}{ C(\eta,m) +\epsilon\lambda} \left(1-e^{-2m\left(\frac 1\epsilon C(\eta,m) +\lambda\right) t}\right)|x_0|^{2m}}.
	\end{align*}
	Denote $\zeta(t)=\frac { C(\eta,m) +\epsilon\lambda}{ (1-\eta)a_{2m+1} \left(1-e^{-2m\left(\frac 1\epsilon C(\eta,m) +\lambda\right)t}\right)}$.
	It follows that for $t>0,$
	\begin{align}\label{x-dep}
		\left|\Phi_t^{\lambda}(x_0)\right|\le \zeta(t)^{\frac 1{2m}}.
	\end{align}
	We complete the proof by combining  \eqref{t-dep} and \eqref{x-dep}.
\fi

\end{proof}

The preceding argument gives the extension to \(\phi\in\mathcal B_{V,q}(\mathbb H)\)
for the strong Feller contribution.  In the asymptotic strong Feller estimate
we keep the assumption \(\phi\in\mathcal C_b^1(\mathbb H)\cap\mathcal B_{V,q}(\mathbb H)\)
because the residual term in \eqref{eq:gradient-final} involves
\(\|D\phi\|_\infty\).  If the full non-degeneracy condition \eqref{erg-noi}
is imposed, this derivative term disappears.

\begin{rk}\label{rk-gamma}
We would like to note that $\gamma(\epsilon)$ typically also depends on $\frac {1} {\epsilon}$ polynomially. 
If, for example, $\lambda_N\gtrsim N^{2/d}$ and $q_i\gtrsim i^{-m_1}$ for some $m_1>0$, then $N_\epsilon\lesssim \epsilon^{-d/2}$ and hence
\[
    \gamma(\epsilon)
    =\sup_{1\le i\le N_\epsilon}q_i^{-1/2}
    \lesssim N_\epsilon^{m_1/2}
    \lesssim \epsilon^{-m_1d/4}.
\]
Thus
\[
    \frac{\gamma(\epsilon)}{\sqrt\epsilon}
    \lesssim \epsilon^{-\frac12-\frac{m_1d}{4}},
\]
which is polynomial in $\epsilon^{-1}$.
\end{rk}

\begin{rk}
There are other sufficient conditions to derive a derivative estimate that depends on $\frac 1\epsilon$ polynomially, for instance, the non-degeneracy condition used in \cite{CHS21}:
\begin{align}\label{erg-noi-chs}
	\rev{\left\|Q^{-\frac 12}(-A+\lambda I)^{-\frac 12}\right\|_{\mathcal L(\mathbb H)}<\infty.}
\end{align}
\rev{This assumption is the non-degeneracy condition used in \cite{CHS21}.  If \(Qe_j=q_je_j\) and \(-Ae_j=\lambda_je_j\), then \eqref{erg-noi-chs} is equivalent to}
\[
\rev{\sup_{j\ge1}q_j^{-1/2}(\lambda_j+\lambda)^{-1/2}<\infty.}
\]
\rev{It is compatible with the Hilbert--Schmidt condition in Theorem \ref{main-tm} for a range of trace-class noises.  For example, if \(q_j\asymp(\lambda_j+\lambda)^{-\rho}\), then \eqref{erg-noi-chs} requires \(\rho\le1\), while \(\|(-A+\lambda)^{(\beta-1)/2}Q^{1/2}\|_{HS(\mathbb H)}<\infty\) requires \(\rho>\beta-1+d/2\).  Hence, admissible examples exist whenever \(\beta<2-d/2\). }

The strong Feller property and irreducible property of $P_t$, together with the time uniform regularity estimate in $\HH^{\beta}$ (i.e., (ii) in Lemma \eqref{exa}), lead to the existence of a unique invariant measure $\mu^*$ for $P_t$. Although this invariant measure is known to be exponentially mixing, i.e., 
\begin{align*}
\sup_{\|\phi\|_{V,q}\le 1} |P_t\phi(x)-\<\phi,\mu^*\>|\le \hat C V(x) e^{- \hat \omega t},
\end{align*}
the constant $\hat C$ or exponent  $\hat \omega$  may still depend on $\frac 1{\epsilon}$ exponentially \cite{GM06}. In particular, by \cite[Theorem 12.5]{GM06}, one cam show that if $q=0,$
then $\omega=-\frac 12\log(1-\widehat \eta)$, $\widehat \eta=\frac 12 c_1 e^{-c_2R^{4m+2}}\int_{\mathbb H}e^{-\Lambda(y)}\mu_1(dy)$,   $R=2(\frac { C(\eta',m) +\epsilon\lambda}{ (1- \eta')a_{2m+1} (1-e^{-2m(\frac 1\epsilon C(\eta,m) +\lambda)\widetilde T})})$ with $C(\eta',m)$ only depending on the time-uniform moment bound of $\|Z(\cdot)\|$ and the power $2m$, and $\widehat C=2(1-\widehat \eta)^{-1},$ with $\widetilde T>0$ large enough and $\eta'\in (0,1)$ small enough.  Here $\Lambda$ is a measurable mapping,  $c_1,c_2$ are the constant coefficients appearing in the lower bound of transition density with respect to a Gaussian measure (for more details, we refer to \cite[Theorem 5.3]{GM06a}). From the proof of  \cite[Theorem 5.3]{GM06a}, it can be seen that \begin{align*}
		c_1&=\exp\left(-C_1\left(\frac {1}{\epsilon}+\lambda\right)^2\right), \; c_2=C_2\left(\frac {1}{\epsilon}+\lambda\right)^2.
	\end{align*}
	As a consequence, one can see that $\widehat C$ is bounded in dependent of $\frac 1{\epsilon}$, but $\widehat \omega \lesssim e^{-\frac c{\epsilon^2}}$ for some $c>0.$ 
\end{rk}

%Choose $\lambda>0$ large enough such that  $\lambda  \sim O(\frac 1{\epsilon})$. 

\subsection{Proof of Theorem \ref{main-kol-reg}}

Now we prove the time-independent Kolmogorov regularity estimate stated in
Theorem \ref{main-kol-reg}.

\begin{proof}[Proof of Theorem \ref{main-kol-reg}]
	\rev{We first treat the short-time case.
	By  Proposition \ref{reg-kl-short}, for \(0<t\le 2T_0\) with \(T_0=2\epsilon\),
	\eqref{reg-kol-1}-\eqref{reg-kol-2} leads to the \revblue{desired} estimate \eqref{time-ind1}--\eqref{time-ind2}. }
	%Thus, no limiting argument is needed to obtain the estimate for fixed \(\delta\); the passage \(\delta\to0\) is used only after the uniform estimate has been established. 
	 %In the short-time estimates below, \(T\) denotes this fixed horizon \(T_0\).
	\rev{It remains to estimate the case \(t>T_0\).} 
	
\textbf{Proof of \eqref{time-ind1}:} 
\iffalse
For $T_0>0$, we denote by $\lfloor t\rfloor_{T_0}$ the largest multiple of $T_0$
that is less than or equal to $t$. Equivalently,
\[
\lfloor t\rfloor_{T_0}=K'T_0,
\qquad
K'=\left\lfloor \frac{t}{T_0}\right\rfloor .
\]
\fi
	According to the Markov property of $P_t$ and \eqref{reg-kol-1}, 
	we have that 
	\begin{align*}
	\|D P_t\phi(x)\|_{\mathbb H^{2\alpha_1}}
	&\le 
	\|DP_{t-\frac {T_0} 2} \phi\|_{\infty} P_{2m(2m+1),2m}(\mathfrak q_\beta,\|x\|_E)\frac 1{\epsilon}\exp\Big(\frac {C'T_0}{\epsilon} \Big)
		(1+T_0^{-\alpha_1}).
\end{align*}
	Furthermore,  applying \eqref{asy-strong-feller} and $T_0=2\epsilon$, we get
	\begin{align}\label{first-der-upper0}
	\|D P_t\phi(x)\|_{\mathbb H^{\alpha_1}}
	&\le P_{2m(2m+1),2m}(\mathfrak q_\beta,\|x\|_E)\frac 1{\epsilon^{1+\alpha_1}}
		\\\nonumber
		&\times 
		[\tilde b(\epsilon^{-1},t-\epsilon, \mathfrak q_\beta,\gamma(\epsilon))\|\phi\|_{V,q}+c_1e^{-\frac {c_0}\epsilon t}\|D\phi\|_{\infty}].
	\end{align}

	\revred{Combining this with \eqref{b-strong-feller} and \eqref{asy-strong-feller1}, and using \(t\ge 2T_0\) and \(\alpha_1\in[0,1)\), we obtain a positive polynomial \(\hat P_{2m(2m+1)+q^2,2m}(\mathfrak q_\beta,\|x\|_E)\) and \(a(\epsilon^{-1})=\gamma(\epsilon)\epsilon^{-5/2}\) such that}
\begin{align}\label{first-oder-large}
&|DP_t\phi(x)\cdot h_1| 
\le \hat P_{2m(2m+1)+q^2,2m}(\mathfrak q_\beta,\|x\|_E) (\|\phi\|_{V,q}+c_1\|D\phi\|_{\infty})a(\frac 1\epsilon)\|h_1\|_{\mathbb H^{-2\alpha_1}}.
\end{align}
This, together with \eqref{reg-kol-1}, leads to \eqref{time-ind1}. 

By similar arguments, one can obtain \eqref{time-ind1} with $c_1=0$ and $a(\frac 1\epsilon)=\mathfrak q_{-}\epsilon^{-5/2}$ under the condition \eqref{erg-noi}.

  \textbf{Proof of \eqref{time-ind2}:}
  Note that $t>2T_0>T_0$. 	
  To prove the second-order derivative estimate, we recall the 
	the Bismut-Elworthy-Li formula (see, e.g., \cite{Cer01}), for $\phi\in  \mathcal B_{V,q}(\mathbb H),$ 
	\begin{align*}
		D P_t \phi(x)\cdot h&=\frac 1 t\mathbb E\Big[\int_0^t \left\<Q^{-\frac 12}\eta^h(s,x),d\widetilde W(s)\right\>\phi(u(t,x))\Big]
	\end{align*}
	 Here $\widetilde W$ is another cylindrical Wiener process, and $h\in \mathbb H.$ 
Thanks to the Markov property of $P_t$, we have that  for $\phi \in \mathcal B_{V,q}(\mathbb H)$ and any $h_1,h_2\in \mathbb H$, 
	\begin{equation}\label{bel-2}
\begin{aligned}
&D^2 P_t\phi(x)\cdot (h_1,h_2)\\
&=\frac{2}{T_0}\mathbb E\Bigg[
\int_0^{\frac {T_0}2}
\left\langle
Q^{-\frac12}\zeta^{h_1,h_2}(s,x),
d\widetilde W(s)
\right\rangle
\,
P_{t-\frac {T_0}2}
\phi(u(\tfrac {T_0}2,x))
\Bigg]
\\
&\quad
+\frac{2}{T_0}\mathbb E\Bigg[
\int_0^{\frac {T_0}2}
\left\langle
Q^{-\frac12}\eta^{h_1}(s,x),
d\widetilde W(s)
\right\rangle
\,
D P_{t-\frac {T_0}2}\phi(u(\tfrac {T_0}2,x))
\cdot
\eta^{h_2}(\tfrac {T_0}2,x)
\Bigg].
\end{aligned}
\end{equation} 
	
\revred{Applying H\"older's inequality to \eqref{bel-2}, we get}
\begin{align*}
&|D^2 P_t\phi(x)\cdot (h_1,h_2)|\\
&\le \frac{2}{T_0}\sqrt{\mathbb E\Bigg[
\int_0^{\frac {T_0}2}
\|Q^{-\frac12}\zeta^{h_1,h_2}(s,x)\|^2ds\Big]}
\sqrt{\E\Big[
|P_{t-\frac {T_0}2}
\phi(u(\tfrac {T_0}2,x))|^2
\Big]}
\\
&\quad
+\frac{2}{T_0}\sqrt{\mathbb E\Bigg[
\int_0^{\frac {T_0}2}
\|Q^{-\frac12}\eta^{h_1}(s,x)\|^2ds }
\,
\E \Big[\Big|D P_{t-\frac {T_0}2}\phi(u(\tfrac {T_0}2,x))
\cdot
\eta^{h_2}(\tfrac {T_0}2,x)\Big|^2\Big]
\Bigg]\\
&=:I_1+I_2.
\end{align*}	

Next, we will bound the two terms $I_1$ and $I_2$, respectively.	

\textbf{Upper bound of $I_1$:}
\revred{Recall that \eqref{erg-noi} implies that $d=1$.}
Notice that by \eqref{var-sec} with $\alpha_1+\alpha_2\in (0,1)$, $\alpha\in (\frac d4,\frac 12)$ and $p=2$, 
\begin{align*}
 &\frac{2}{T_0}\sqrt{\mathbb E\Bigg[
\int_0^{\frac {T_0}2}
\|Q^{-\frac12}\zeta^{h_1,h_2}(s,x)\|^2ds\Big]}\\
&\le C \mathfrak q_{-} \frac{2}{T_0} \left(\frac1\epsilon\right)^{4} \left(1+\|x\|_E^{8m-1}\right) \sqrt{\mathbb E\Bigg[
\int_0^{\frac {T_0}2}
		(1+s^{2(1-\alpha-\alpha_1-\alpha_2)})ds\Big]} P_{(8m-1)(2m+1)}(\mathfrak q_{\beta})\\
		&\quad \times \|h_1\|_{\mathbb H^{-\alpha_1}}\|h_2\|_{\mathbb H^{-\alpha_2}}\\
		&\le C\mathfrak q_{-} P_{(8m-1)(2m+1)}(\mathfrak q_{\beta}) T_0^{-1}\left(\frac1\epsilon\right)^{4} \left(1+\|x\|_E^{8m-1}\right) \|h_1\|_{\mathbb H^{-\alpha_1}}\|h_2\|_{\mathbb H^{-\alpha_2}}\\
		&\le C\mathfrak q_{-} P_{(8m-1)(2m+1)}(\mathfrak q_{\beta}) \left(\frac1\epsilon\right)^{5} \left(1+\|x\|_E^{8m-1}\right) \|h_1\|_{\mathbb H^{-\alpha_1}}\|h_2\|_{\mathbb H^{-\alpha_2}}.
\end{align*}

On the other hand, by \eqref{eq:initial-free-E} and the fact that $t>T_0$
\begin{align*}
\sqrt{\E\Big[
|P_{t-\frac {T_0}2}
\phi(u(\tfrac {T_0}2,x))|^2
\Big]}
&\le C \|\phi\|_{V,q} P_{q}(\mathfrak q_{\beta})^{q}\sqrt{\E\Big[(1+
\|u(t,x)\|^{2q})
\Big]}\\
&\le C \|\phi\|_{V,q} P_{q}(\mathfrak q_{\beta})^{q} \left(1+\left(\frac{\epsilon}{T_0}\right)^{q/(2m)}\right)\le C \|\phi\|_{V,q} P_{q^2}(\mathfrak q_{\beta}).
\end{align*}

Combining the above two estimates, we get 
\begin{align*}
I_1\le C   \mathfrak q_{-} P_{(8m-1)(2m+1)+q^2}(\mathfrak q_{\beta}) \|\phi\|_{V,q} \left(\frac1\epsilon\right)^{5} \left(1+\|x\|_E^{8m-1}\right)  \|h_1\|_{\mathbb H^{-\alpha_1}}\|h_2\|_{\mathbb H^{-\alpha_2}}.
\end{align*}

\textbf{Upper bound of $I_2$:}	
	By \eqref{smo-two} with $\alpha_1=0$, we have 
	 
\begin{align*}
	&\frac{2}{T_0}\sqrt{\mathbb E\Big[
\int_0^{\frac {T_0}2}
\|Q^{-\frac12}\eta^{h_1}(s,x)\|^2ds \Big]}\\&\le C\mathfrak q_{-}\frac{2}{T_0}\frac1\epsilon \left(1+\|x\|_E^{2m}\right)P_{2m(2m+1)}(\mathfrak q_{\beta})  \sqrt{\mathbb E\Big[
\int_0^{\frac {T_0}2}
\left(1+e^{-c_0s}\right) ds \Big]} \|h_1\|\\
&\le  C\mathfrak q_{-}\frac1 {\epsilon^{\frac 32}} P_{2m(2m+1)}(\mathfrak q_{\beta})  \left(1+\|x\|_E^{2m}\right) \|h_1\|.
\end{align*}

On the other hand, using \eqref{asy-strong-feller} with $c_1=0,$ and \eqref{smo-two} again, we get that for any $\alpha_2\in [0,1),$
\begin{align*}
&\sqrt{\E \Big[\Big|D P_{t-\frac {T_0}2}\phi(u(\tfrac {T_0}2,x))
\cdot
\eta^{h_2}(\tfrac {T_0}2,x)\Big|^2\Big]}\\
&\le  \tilde b(\frac 1\epsilon,{t-\frac {T_0}2},\mathfrak q_{\beta},\mathfrak q_{-}) \|\phi\|_{V,q} \sqrt{\E \|\eta^{h_2}(\tfrac {T_0}2,x)\|^2}\\
&\le C \tilde b(\frac 1\epsilon,{t-\frac {T_0}2},\mathfrak q_{\beta},\mathfrak q_{-})  \|\phi\|_{V,q} \frac1 {\epsilon^{1+\alpha_2}}\left(1+\|x\|_E^{2m}\right) P_{2m(2m+1)}(\mathfrak q_{\beta})       \|h_2\|_{\mathbb H^{-2\alpha_2}}.
\end{align*}

Thus, we conclude that 
\begin{align*}
I_2&\le  C \tilde b(\frac 1\epsilon,{t-\frac {T_0}2},\mathfrak q_{\beta},\mathfrak q_{-}) \|\phi\|_{V,q} \frac1 {\epsilon^{\frac 52+\alpha_2}}P_{4m(2m+1)}(\mathfrak q_{\beta}) \left(1+\|x\|_E^{4m}\right)  \|h_1\|\|h_2\|_{\mathbb H^{-2\alpha_2}}.
\end{align*}	

Combining the estimates of $I_1$ with $\alpha_1=\alpha_2=0$ and $I_2$ with $\alpha_2=0$, we obtain 	
	\begin{align}\label{second-der-upper0}
	&|D^2 P_t\phi(x)\cdot (h_1,h_2)|\\\nonumber
	&\le  C P_{(8m-1)(2m+1)}(\mathfrak q_{\beta})  \tilde b(\frac 1\epsilon,{t-\frac {T_0}2},\mathfrak q_{\beta},\mathfrak q_{-}) \|\phi\|_{V,q} \frac 1{\epsilon^{5}} \left(1+\|x\|_E^{8m-1}\right) \|h_1\|\|h_2\|.
	\end{align}

Next, we prove the improved regularity estimate for the second-order derivative.  
By the Markov property, we further get that for $t>2T_0,$
\begin{align}\label{second-der-upper1}
D^2P_t\phi(x)\cdot(h_1,h_2)
&=
\mathbb E\left[
D^2P_{t-\frac {T_0}2}\phi\big(u(\frac {T_0}2,x)\big)
\cdot
\big(
\eta^{h_1}(\frac {T_0}2,x),
\eta^{h_2}(\frac {T_0}2,x)
\big)
\right]
\\\nonumber
&\quad+
\mathbb E\left[
DP_{t-\frac {T_0}2}\phi\big(u(\frac {T_0}2,x)\big)
\cdot
\zeta^{h_1,h_2}(\frac {T_0}2,x)
\right].
\end{align}
This, together with \eqref{second-der-upper0}, \eqref{first-oder-large} and \eqref{eq:initial-free-E0}, as well as  \eqref{var-sec} and \eqref{smo-two}, yields that 
\begin{align*}&|D^2P_t\phi(x)\cdot(h_1,h_2)|\\
&\le C \|\phi\|_{V,q}P_{(8m-1)(2m+1)}(\mathfrak q_{\beta})  \tilde b(\frac 1\epsilon,{t-\frac {T_0}2},\mathfrak q_{\beta},\mathfrak q_{-}) \frac 1{\epsilon^{5}} \revred{\sqrt{\E \left[1+\|u(\frac {T_0}2,x)\|_E^{8m-1}\right]}}
		\\
		&\quad \times \|\eta^{h_1}(\frac {T_0}2,x)\|_{L^4(\Omega;\mathbb H)}\|\eta^{h_2}(\frac {T_0}2,x)\|_{L^4(\Omega;\mathbb H)}\\
&+ P_{4m(2m+1)+q^2,2m}(\mathfrak q_\beta,\|x\|_E) \|\phi\|_{V,q}a(\frac 1\epsilon)\|\zeta^{h_1,h_2}(\frac {T_0}2,x)\|_{L^2(\Omega;\mathbb H)}\\
		&\le C\|\phi\|_{V,q} P_{(8m-1)(2m+1)}(\mathfrak q_{\beta})  \tilde b(\frac 1\epsilon,{t-\frac {T_0}2},\mathfrak q_{\beta},\mathfrak q_{-})  \frac 1{\epsilon^{5}} \sqrt{\E(1+\|u(\frac {T_0}2,x)\|_E^{8m-1})}	\\
		&\quad \times \frac1 {\epsilon^2}P_{4m(2m+1)}(\mathfrak q_{\beta}) \left(1+\|x\|_E^{4m}\right)  \|h_1\|_{\mathbb H^{-2\alpha_1}}\|h_2\|_{\mathbb H^{-2\alpha_2}}\left(1+T_0^{-\alpha_1-\alpha_2}\right).\\
&+ P_{4m(2m+1)+q^2,2m}(\mathfrak q_\beta,\|x\|_E)  
 \|\phi\|_{V,q}\frac 1 {\epsilon^4} a(\frac 1\epsilon) P_{(8m-1)(2m+1),8m-1}(\mathfrak q_\beta,\|x\|_E)   \\
		&\quad \times (1+T_0^{1-\alpha-\alpha_1-\alpha_2})\|h_1\|_{\mathbb H^{-2\alpha_1}}\|h_2\|_{\mathbb H^{-2\alpha_2}}.
\end{align*}

By  \eqref{eq:initial-free-E} with the time being $T_0$, one further has that 
\begin{align*}
&|D^2P_t\phi(x)\cdot(h_1,h_2)|\\
&\le 
C\|\phi\|_{V,q}\tilde b(\frac 1\epsilon,{t-\frac {T_0}2},\mathfrak q_{\beta},\mathfrak q_{-})  \frac 1{\epsilon^{8}} P_{(12m-1)(2m+1)+4m(8m-1)}(\mathfrak q_\beta)\left(1+\|x\|_E^{4m}\right)\\
&\quad \times	  \|h_1\|_{\mathbb H^{-2\alpha_1}}\|h_2\|_{\mathbb H^{-2\alpha_2}}\\
&+ P_{(2m+1)(10m-1)+q^2,10m-1}(\mathfrak q_\beta,\|x\|_E)  
 \|\phi\|_{V,q}\frac 1{\epsilon^4} a(\frac 1\epsilon) 
		\|h_1\|_{\mathbb H^{-2\alpha_1}}\|h_2\|_{\mathbb H^{-2\alpha_2}}.
\end{align*}

This, together with \eqref{b-strong-feller} and \eqref{asy-strong-feller1}, as well as the condition that $t\ge T_0,$ we conclude that 
\begin{align*}
|D^2P_t\phi(x)\cdot(h_1,h_2)| 
&\le \mathfrak q_{-}\|\phi\|_{V,q} \frac {1}{\epsilon^{\frac {17}2}} P_{(12m-1)(2m+1)+4m(8m-1)+q^2,10m-1}(\mathfrak q_\beta,\|x\|_E) \\
&\quad \times \|h_1\|_{\mathbb H^{-2\alpha_1}}\|h_2\|_{\mathbb H^{-2\alpha_2}},
\end{align*}
which, together with \eqref{reg-kol-2}, yields  \eqref{time-ind2}.
\end{proof}

\section{Proof of Theorem \ref{main-tm}}
\label{sec-5}

\revblue{This section proves Theorem \ref{main-tm}.  The proof shows how the time-independent Kolmogorov regularity estimates replace the deterministic spectral estimate used in sharp-interface analyses of the deterministic Allen--Cahn equation.  We carry out the argument for the splitting approximation \eqref{spl-sch}; the same structure can be used for other discretizations once the corresponding time-uniform moment and regularity estimates are available.}

\revblue{The proof has two steps.  First, we introduce a continuous interpolation of the splitting scheme, following the standard weak-error strategy in \cite{BG18B,CHS21}.  Second, we insert this interpolation into the Kolmogorov equation and estimate the resulting local defects by the regularity estimates of Section \ref{sec-4}.}

\subsection{Continuous interpolation of splitting scheme}

In this part, we define the continuous interpolation of \eqref{spl-sch} and present its several useful properties.  For any $k\in \mathbb N$ and $t\in[t_k,t_{k+1}],$  the continuous interpolation of \eqref{spl-sch} is defined by
\begin{align*}
	\widetilde u(t)=S_{\lambda}(t-t_k)\Phi_{t-t_k}^{\lambda}(u_{k})+\int_{t_k}^tS_{\lambda}(t-s)dW(s),
\end{align*} 
which satisfies 
\begin{align*}
	d \widetilde u(t)=(A-\lambda)\widetilde u dt+S_{\lambda}(t-t_k)\Psi_{0}^{\lambda}(\Phi_{t-t_k}^{\lambda}(u_k))dt+dW(t),\; t\in[t_k,t_{k+1}].
\end{align*}
\revblue{The same argument also applies to non-uniform grids, provided that the mesh sizes are uniformly bounded by a constant of order \(\epsilon\).}

\revred{By the a priori estimate for \eqref{spl-sch} (see Lemma \ref{pri-un-e}), the continuous interpolation \(\widetilde u(t)\) satisfies the following time-independent bound: for any \(p>0\),}
\begin{align}\label{pri-con}
	\sup_{t\ge 0}\left\|\widetilde u(t)\right\|_{L^p(\Omega;E)}\le \rev{P_{2m+1}(\mathfrak q_\beta)}(1+\|u(0)\|_E).
\end{align}

\iffalse
\revblue{For later use we also record the elementary derivative bounds for
the original drift}
\[
\revblue{\Theta_0^\lambda(\xi):=\Psi_0^\lambda(\xi)
=-\epsilon^{-1}f(\xi)+\lambda \xi .}
\]
\begin{align}\label{theta0-der}
\revblue{|(\Theta_0^\lambda)'(\xi)|\le C(\epsilon^{-1}+\lambda)(1+|\xi|^{2m}),\qquad
|(\Theta_0^\lambda)''(\xi)|\le C\epsilon^{-1}(1+|\xi|^{2m-1}).}
\end{align}

\iffalse
\rev{For the original Kolmogorov equation, we use the drift}
\[
\rev{\Psi_0^\lambda(\xi)=-\epsilon^{-1}f(\xi)+\lambda\xi.}
\]
\rev{Its derivatives satisfy}
\begin{align}\label{theta0-der}
\rev{|(\Theta_0^\lambda)'(\xi)|\le C(\epsilon^{-1}+\lambda)(1+|\xi|^{2m}),\qquad
|(\Theta_0^\lambda)''(\xi)|\le C\epsilon^{-1}(1+|\xi|^{2m-1}).}
\end{align}
\fi

\fi

\revred{We also need the following negative Sobolev multiplier estimate to obtain the sharp weak convergence rate with respect to \(\tau\). Its proof is given in Appendix \ref{negative-multiplier}.}

\begin{lm}[A negative Sobolev multiplier estimate]\label{lem:negative-multiplier}
Let \(s\in(\frac d2,2)\) and $s\ge \gamma_1>0$.  For \(x\in E\cap\mathbb H^{\gamma_1}\) and
\(y\in\mathbb H^{-\gamma_1}\) , there exists a constant \(C>0\) such that
\[
 \|\revred{(-\frac 1\epsilon f'(x)+\lambda)}y\|_{\mathbb H^{-s}}
 \le
 C\epsilon^{-1}
 \left(1+\|x\|_E^{2m}+\|x\|_{\mathbb H^{\gamma_1}}^{2m}\right)
 \|y\|_{\mathbb H^{-\gamma_1}} .
\]
\end{lm}

\subsection{Polynomial dependence in weak error analysis}
\;
Let \(X(t,x)=P_t\phi(x)=\E[\phi(u(t,x))]\) denote the solution of the
Kolmogorov equation.  
\rev{Taking \(T=N\tau\), the Kolmogorov equation and the interpolation
\(\widetilde u\) give}
\begin{align*}
	&\rev{\E[\phi(u(T,x))]-\E[\phi(u_N)]}\\
	&=\sum_{k=0}^{N-1} \left(\E\left[X(T-t_k, u_k)\right]-\E\left[X(T-t_{k+1},u_{k+1})\right]\right)\\
	&=\sum_{k=0}^{N-1} \int_{t_k}^{t_{k+1}} \E \Big[\left\< DX(T-t, \widetilde u(t)),
	\Psi_{0}^{\lambda}(\tilde u(t))-S_{\lambda}(t-t_k)\Psi_0^{\lambda}\left(\Phi^{\lambda}_{t-t_k}(u_k)\right)\right\> \Big]dt.
\end{align*}

\begin{proof}[Proof of Theorem \ref{main-tm}]
Our strategy is to decompose the above weak error as 
\begin{align*}
	&\sum_{k=0}^{N-1} \int_{t_k}^{t_{k+1}} \E \Big[ \left\< DX(T-t, \widetilde u(t)), \Psi_{0}^{\lambda}(\tilde u(t))-S_{\lambda}(t-t_k)\Psi_0^{\lambda}\left(\Phi^\lambda_{t-t_k}(u_k)\right)\right\> \Big]dt\\
	&=\sum_{k=0}^{N-1} \Big(\int_{t_k}^{t_{k+1}} \E \Big[\left\< DX(T-t, \widetilde u(t)),\Psi_{0}^{\lambda}(\tilde u(t))-\Psi_{0}^{\lambda}(u_k)\right\>\Big]dt\\
	&+\int_{t_k}^{t_{k+1}} \E \Big[\left\< DX(T-t, \widetilde u(t)), \left( I-S_{\lambda}(t-t_k)\right)\Psi_0^{\lambda}(u_k) \right\> \Big]dt\\
	&+\int_{t_k}^{t_{k+1}} \E \Big[\left\< DX(T-t, \widetilde u(t)), S_{\lambda}(t-t_k)\left(\Psi_0^{\lambda}(u_k)-\Psi_0^{\lambda}\left(\Phi^{\lambda}_{t-t_k}(u_k)\right)\right) \right\> \Big]dt\Big)\\
	&=:\sum_{k=0}^{N-1} Err_1^{k+1}+Err_2^{k+1}+Err_3^{k+1}.
\end{align*}
\rev{Next, we estimate \(Err_1^{k+1}\)-\(Err_3^{k+1}\) using the time-independent
regularity estimate in Theorem \ref{main-kol-reg}.}

\textbf{Estimate of $Err_1^{k+1}$:}
%:
We decompose $Err_1^{k+1}$ by 

\begin{align*}
	|Err_1^{k+1}|
	&\le \Big|\int_{t_k}^{t_{k+1}} \E \Big[ \left\< DX(T-t, \widetilde u(t)), \revred{(-\frac 1\epsilon f'(u_k)+\lambda)}\left(\widetilde u(t)-u_k\right)\right\>\Big]dt\Big|\\
	&+ \Big|\int_{t_k}^{t_{k+1}} \E \Big[ \big\< DX(T-t, \widetilde u(t)), \revred{-}\frac 1\epsilon \int_0^1 f''\left((1-\theta)\widetilde u(t)+\theta u_k\right)d\theta \cdot
	\\
	&\quad \left(\widetilde u(t)-u_k,\widetilde u(t)-u_k \right)\big\>\Big]dt\Big|\\
	&=:Err_{11}^{k+1}+Err_{12}^{k+1}.
\end{align*}
For $Err_{11}^{k+1}$, the martingale property of the stochastic integral and the adaptivity of $\widetilde u$ yield that
\begin{align*}
	&|Err_{11}^{k+1}|\\
	&\le \Big|\int_{t_k}^{t_{k+1}} \E \Big[ \left\< DX(T-t, \widetilde u(t)), \revred{(-\frac 1\epsilon f'(u_k)+\lambda)}\cdot \left(S_{\lambda}(t-t_k)-I\right)u_k\right\>\Big]dt\Big|\\
	&+\Big|\int_{t_k}^{t_{k+1}} \E \Big[ \left\< DX(T-t, \widetilde u(t)), \revred{(-\frac 1\epsilon f'(u_k) +\lambda )}\cdot S_{\lambda}(t-t_k)(\Phi_{t-t_k}^{\lambda}(u_k)-u_k)\right\>\Big] dt\Big|\\
	&=: Err_{111}^{k+1}+Err_{112}^{k+1}.
\end{align*}

{Choose \(\gamma\in(\max\{d/4,\beta/2\},1)\), which is possible since \(d\le3\) and \(\beta<2\).  Applying \eqref{time-ind1} to \(Err_{111}^{k+1}\) with \(\alpha_1=\gamma\), and using Lemma \ref{lem:negative-multiplier} with \(s=2\gamma\) and \(\gamma_1=\beta\), gives the desired negative-Sobolev estimate.}
Together with Lemma \ref{pri-un-e}, H\"older's inequality, and \eqref{eq:initial-free-E1}, we \revred{obtain}
\begin{align*}
&|Err_{111}^{k+1}|\\
&\le\int_{t_k}^{t_{k+1}}    P_{2m(2m+1)+q^2,2m}(\mathfrak q_\beta,\|u(0)\|_E) (\|\phi\|_{V,q}+c_1\|D\phi\|_{\infty})a(\frac 1\epsilon) \\
&\quad  \times \Big([1+(T-t)^{-\gamma}]\left\|\revred{(-\frac 1\epsilon f'(u_k)+\lambda)}\cdot \rev{\left(S_\lambda(t-t_k)-I\right)u_k}\right\|_{L^2(\Omega;\mathbb H^{-2\gamma})} dt\\
&\le \int_{t_k}^{t_{k+1}}    P_{2m(2m+1)+q^2,2m}(\mathfrak q_\beta,\|u(0)\|_E)  (\|\phi\|_{V,q}+c_1\|D\phi\|_{\infty})a(\frac 1\epsilon) \frac 1{\epsilon^{1+2m}} [1+(T-t)^{-\gamma}] \\
	&\quad \times \rev{P_{(2m+1)^2,2m+1,1}(\mathfrak q_\beta, \|u(0)\|_E,\|u(0)\|_{\mathbb H^{\beta}})}^{2m}\left\|\rev{\left(S_\lambda(t-t_k)-I\right)u_k}\right\|_{L^4(\Omega;\mathbb H^{-\beta})} dt.
\end{align*}
It follows from \eqref{smooth0} and \eqref{reg-nu} that
\begin{align*}
&|Err_{111}^{k+1}|\\
&\le \int_{t_k}^{t_{k+1}}    \rev{P_{2m(2m+1)(2m+2)+q^2,2m(2m+2),2m}(\mathfrak q_\beta, \|u(0)\|_E,\|u(0)\|_{\mathbb H^{\beta}})}(\|\phi\|_{V,q}+c_1\|D\phi\|_{\infty})\\
	&\quad \times a(\frac 1\epsilon) \frac 1{\epsilon^{1+2m}} [1+(T-t)^{-\gamma}] \tau^{\beta} \left\|u_k\right\|_{L^4(\Omega;\mathbb H^{\beta})} dt\\
&\le \int_{t_k}^{t_{k+1}}    \rev{P_{(2m+1)(4m^2+6m+1)+q^2,2m(2m+2)+2m+1,2m+1}(\mathfrak q_\beta, \|u(0)\|_E,\|u(0)\|_{\mathbb H^{\beta}})}\\
	&\quad \times  a(\frac 1\epsilon) \frac 1{\epsilon^{2+2m}}  (\|\phi\|_{V,q}+c_1\|D\phi\|_{\infty}) [1+(T-t)^{-\gamma}]   \tau^{\beta }dt.
\end{align*}

Then applying \eqref{time-ind1} to $Err_{112}^{k+1}$ with $\alpha_1=0,$ by similar arguments, one can obtain  that
\begin{align*}
|Err_{112}^{k+1}|&\le\int_{t_k}^{t_{k+1}}    P_{2m(2m+1)+q^2,2m}(\mathfrak q_\beta,\|u(0)\|_E) (\|\phi\|_{V,q}+c_1\|D\phi\|_{\infty})a(\frac 1\epsilon)\\
	&\quad \times 	\Big\|\revred{(-\frac 1\epsilon f' (u_k) +\lambda)}\cdot \rev{S_\lambda(t-t_k)\left(\Phi_{t-t_k}^{\lambda}(u_k)-u_k\right)}\Big\|_{L^2(\Omega;\mathbb H)}dt\\
	&\le \int_{t_k}^{t_{k+1}}     P_{4m(2m+1)+q^2,4m}(\mathfrak q_\beta,\|u(0)\|_E)  (\|\phi\|_{V,q}+c_1\|D\phi\|_{\infty})a(\frac 1\epsilon)\frac 1\epsilon \\
	&\quad \times\Big\|\Phi_{t-t_k}^{\lambda}(u_k)-u_k\Big\|_{L^4(\Omega;\mathbb H)}dt.
\end{align*}

According to Lemma \ref{pri-un-e}  and the property \eqref{prop-psi} of $\Psi_{t}^{\lambda}$, thanks to $\tau \lesssim \epsilon$, we obtain that
\begin{align*}
&|Err_{112}^{k+1}|\\
 &\le \int_{t_k}^{t_{k+1}}     P_{4m(2m+1)+q^2,4m}(\mathfrak q_\beta,\|u(0)\|_E)  (\|\phi\|_{V,q}+c_1\|D\phi\|_{\infty})a(\frac 1\epsilon)\frac 1{\epsilon^2} \\
	&\quad \times \tau \|u_k\|_{L^{8m+4}(\Omega;E)}^{2m+1}dt\\
	&\le \int_{t_k}^{t_{k+1}}     P_{(6m+1)(2m+1)+q^2,6m+1}(\mathfrak q_\beta,\|u(0)\|_E)  (\|\phi\|_{V,q}+c_1\|D\phi\|_{\infty})a(\frac 1\epsilon)\frac 1{\epsilon^2} dt.
\end{align*}

Combining the estimates of $Err_{111}^{k+1}$ and $Err_{112}^{k+1},$ we get 
\begin{align}\nonumber
|Err_{11}^{k+1}|&\le \int_{t_k}^{t_{k+1}}    \rev{P_{(2m+1)(4m^2+6m+1)+q^2,2m(2m+2)+2m+1,2m+1}(\mathfrak q_\beta, \|u(0)\|_E,\|u(0)\|_{\mathbb H^{\beta}})}\\\label{reg-err11}
	&\quad \times  a(\frac 1\epsilon) \frac 1{\epsilon^{2+2m}}  (\|\phi\|_{V,q}+c_1\|D\phi\|_{\infty}) [1+(T-t)^{-\gamma}]  \tau^{ \beta }dt.
\end{align}

Next, we deal with $Err_{12}^{k+1}.$
Applying \eqref{time-ind1} with $\alpha_1\in (\frac d4,1)$ and Lemma \ref{pri-un-e}, together with the dual formulation of the Sobolev embedding theorem and H\"older's inequality, leads to 
\begin{align}\nonumber
	&|Err_{12}^{k+1}|\\\nonumber
	&\le  \Big|\int_{t_k}^{t_{k+1}} \E \Big[ \big\< DX(T-t, \widetilde u(t)), \int_0^1 \frac 1\epsilon f''\left((1-\theta)\widetilde u(t)+\theta u_k\right)d\theta  \\\nonumber
	&\quad \cdot \left(\widetilde u(t)-u_k,\widetilde u(t)-u_k \right)\big\>\Big]dt\Big|\\\nonumber
	&\le  P_{2m(2m+1)+q^2,2m}(\mathfrak q_\beta,\|u(0)\|_E)  \int_{t_k}^{t_{k+1}} (1+(T-t)^{-\alpha_1})(\|\phi\|_{V,q}+c_1\|D\phi\|_{\infty})a(\frac 1\epsilon) \\\nonumber
	&\quad \times \frac 1\epsilon P_{(2m-1)(2m+1),2m-1}(\mathfrak q_\beta,\|u(0)\|_E)\|\widetilde u(t)-u_k \|_{L^8(\Omega;\mathbb H)}^2 dt\\\nonumber
	&=P_{(4m-1)(2m+1)+q^2,4m-1}(\mathfrak q_\beta,\|u(0)\|_E)  \int_{t_k}^{t_{k+1}} (1+(T-t)^{-\alpha_1})(\|\phi\|_{V,q}+c_1\|D\phi\|_{\infty})\\\label{continuity-int}
	&\quad \times a(\frac 1\epsilon) \frac 1\epsilon\|\widetilde u(t)-u_k \|_{L^8(\Omega;\mathbb H)}^2 dt.
\end{align}

Similar to the proof of \eqref{3.8}, one can obtain that for any $p\ge 1,$
\begin{align*}
&\|\widetilde u(t)-u_k \|_{L^p(\Omega; \mathbb H)}\\
&\le \|S_{\lambda}(t-t_k)\Phi_{t-t_k}^{\lambda}(u_{k})-u_k \|_{L^p(\Omega; \mathbb H)}+\|\int_{t_k}^tS_{\lambda}(t-s)dW(s)\|_{L^p(\Omega; \mathbb H)}\\
&\le   \|(S_{\lambda}(t-t_k)-I)u_k \|_{L^p(\Omega; \mathbb H)} +
 \|S_{\lambda}(t-t_k)(\Phi_{t-t_k}^{\lambda}(u_{k})-u_k)\|_{L^p(\Omega; \mathbb H)} \\
 &\quad+C \mathfrak q_{\beta} \tau^{\min (\frac \beta 2,\frac 12)}.
\end{align*}
By \eqref{smooth1} and \eqref{prop-psi}, as well as \eqref{reg-nu} and Lemma \ref{pri-un-e},  it holds that 
\begin{align}\nonumber
&\|\widetilde u(t)-u_k \|_{L^p(\Omega; \mathbb H)}\\\nonumber
&\le   C \tau^{\frac \beta 2}\|u_k \|_{L^p(\Omega; \mathbb H^{\beta})} +
 C\tau \frac 1\epsilon \|u_k\|_{L^{(2m+1)p}(\Omega;E)}^{2m+1}+C \mathfrak q_{\beta} \tau^{\min (\frac \beta 2,\frac 12)}\\\label{continuity-int1}
 &\le \frac 1\epsilon P_{(2m+1)^2, 2m+1,1}(\mathfrak q_{\beta},\|u(0)\|_E,\|u(0)\|_{\mathbb H^{\beta}})\tau^{\min(\frac \beta 2,\frac 12)}.
\end{align}
Substituting \eqref{continuity-int1} with $p=8$ into \eqref{continuity-int}, we   conclude that  
\begin{align}\nonumber
	|Err_{12}^{k+1}|
	&\le P_{(8m+1)(2m+1)+q^2,8m+1,2}(\mathfrak q_\beta,\|u(0)\|_E,\|u(0)\|_{\mathbb H^{\beta}}) a(\frac 1\epsilon) \frac 1{\epsilon^3} \\\label{reg-err-fir}
	&\quad\times (\|\phi\|_{V,q}+c_1\|D\phi\|_{\infty}) \int_{t_k}^{t_{k+1}} (1+(T-t)^{-\alpha_1})   dt \tau^{\min(1,\beta)}.
\end{align}

Combining the estimates \eqref{reg-err-fir} and \eqref{reg-err11} together, we conclude that 

\begin{align}\label{reg-err-fir1}	 &\rev{|Err_{1}^{k+1}|}\\\nonumber
	&\le \int_{t_k}^{t_{k+1}}    \rev{P_{(2m+1)(4m^2+6m+1)+q^2,2m(2m+2)+2m+1,2m+1}(\mathfrak q_\beta, \|u(0)\|_E,\|u(0)\|_{\mathbb H^{\beta}})}\\\nonumber
	&\quad \times  a(\frac 1\epsilon) \frac 1{\epsilon^{2+2m}}  (\|\phi\|_{V,q}+c_1\|D\phi\|_{\infty}) [1+(T-t)^{-\gamma}+(T-t)^{-\alpha_1}] dt  \tau^{\min(\beta,1)}.
	\end{align}

\textbf{Estimate of $Err_2^{k+1}$:}
Applying \eqref{time-ind1} with the parameter $\alpha_3 \in (0,1)$ and \eqref{smooth1}, and using the a priori bound of $u_k$ in Lemma \ref{pri-un-e}, we have
\begin{align*}
	&\int_{t_k}^{t_{k+1}} \E \Big[\left\< DX(T-t, \widetilde u(t)),( I-S_{\lambda}(t-t_k))[\frac 1\epsilon f(u_k)+\lambda u_k]  \right\> \Big]dt\\
	&\le P_{2m(2m+1)+q^2,2m}(\mathfrak q_{\beta},\|u(0)\|_E)\int_{t_k}^{t_{k+1}} (1+(T-t)^{-\alpha_3})(\|\phi\|_{V,q}+c_1\|D\phi\|_{\infty}) a(\frac 1\epsilon) \\\nonumber 
	&\quad   \Big\| \left( I-S_{\lambda}(t-t_k)\right)[\frac 1\epsilon f(u_k)+\lambda u_k]\Big\|_{L^2(\Omega;\mathbb H^{-2\alpha_3})}dt\\
	&\le P_{(4m+1)(2m+1)+q^2,4m+1}(\mathfrak q_{\beta},\|u(0)\|_E)\int_{t_k}^{t_{k+1}} \left(1+(T-t)^{-\alpha_3}\right)dt a(\frac 1\epsilon)\frac 1{\epsilon} 
	\tau^{\alpha_3}.
\end{align*}

\textbf{Estimate of $Err_3^{k+1}$:}
Exploiting \eqref{time-ind1} with $\alpha_1=0$, H\"older's inequality,
Lemma \ref{pri-un-e}, and \eqref{prop-psi},  we have
\begin{align*}
	|Err_{3}^{k+1}|&\le \Big|\int_{t_k}^{t_{k+1}} \E \Big[\left\< DX(T-t, \widetilde u(t)), S_{\lambda}(t-t_k)\left(\Psi_{0}^{\lambda}(u_k)-\Psi_{0}^{\lambda}\left(\Phi_{t-t_k}^{\lambda}(u_k)\right)\right) \right\> \Big]dt\Big|\\
	&\le P_{2m(2m+1)+q^2,2m}(\mathfrak q_{\beta},\|u(0)\|_E) (\|\phi\|_{V,q}+c_1\|D\phi\|_{\infty}) a(\frac 1\epsilon)  \\\nonumber 
	&\quad\times \int_{t_k}^{t_{k+1}}
	\left\|\Psi_{0}^{\lambda}(u_k)-\Psi_{0}^{\lambda}(\Phi_{t-t_k}^{\lambda}(u_k))\right\|_{L^2(\Omega;\mathbb H)}dt\\
	&\le P_{2m(2m+1)+q^2,2m}(\mathfrak q_{\beta},\|u(0)\|_E)  \tau(\|\phi\|_{V,q}+c_1\|D\phi\|_{\infty}) a(\frac 1\epsilon)\frac {1}{\epsilon}  \\
	&\quad \times (1+ \|u_k\|_{L^{8m}(\Omega; E)}^{2m}) \|\Phi_{t-t_k}^{\lambda}(u_k)-u_k\|_{L^4(\Omega;H)}\tau\\
	&\le P_{4m(2m+1)+q^2,4m}(\mathfrak q_{\beta},\|u(0)\|_E)  \tau(\|\phi\|_{V,q}+c_1\|D\phi\|_{\infty}) a(\frac 1\epsilon)\frac {1}{\epsilon}  \\
	&\quad \times \|\Phi_{t-t_k}^{\lambda}(u_k)-u_k\|_{L^4(\Omega;H)}\tau.
\end{align*}
\revred{Together with Lemma \ref{pri-un-e} and the estimate \eqref{prop-psi} for \(\Psi_t^\lambda\), this gives}
\begin{align*}
	|Err_{3}^{k+1}|
	&\le P_{(6m+1)(2m+1)+q^2,6m+1}(\mathfrak q_{\beta},\|u(0)\|_E)  \tau(\|\phi\|_{V,q}+c_1\|D\phi\|_{\infty}) a(\frac 1\epsilon)\frac {1}{\epsilon^2}   \tau^2.
\end{align*}

%\questionpurple{Please check the powers of \(\tau\) in this estimate of \(Err_3^{k+1}\).  The display above appears to contain one factor from the time integral and another factor after \(\|\Phi_{t-t_k}^{\lambda}(u_k)-u_k\|\), before using the flow consistency estimate.  If both factors are intentional, this term is higher order and harmless for the final rate; otherwise one \(\tau\) may be extra.}

Combining the above estimates on $Err_1^{k+1}$-$Err_3^{k+1}$, summing over
$k=0,\ldots, N-1$, and using the integrability of the singular kernel
\((T-t)^{-\alpha'}\) for any \(\alpha' \in(0,1)\), we obtain
\[
	\left|\E\left[\phi(u(t_N))\right]-\E\left[\phi(u_N)\right]\right|
	\le
	(1+T)b_1(\|u(0)\|_E,\|u(0)\|_{\mathbb H^{\beta}}, \epsilon^{-1},\mathfrak q_\beta,\gamma(\epsilon))
	\tau^{\min\left(\beta,\alpha_3\right)} ,
\]
where 
\begin{align}\label{form-poly}
&b_1(\|u(0)\|_E,\|u(0)\|_{\mathbb H^{\beta}}, \epsilon^{-1},\mathfrak q_\beta,\gamma(\epsilon))\\\nonumber 
&= \rev{P_{(2m+1)(4m^2+6m+1)+q^2,2m(2m+2)+2m+1,2m+1 }(\mathfrak q_\beta, \|u(0)\|_E,\|u(0)\|_{\mathbb H^{\beta}})} a(\frac 1\epsilon) \frac 1{\epsilon^{2+2m}} \end{align}
with $a(\frac 1\epsilon)=\frac 1{\epsilon^{\frac 52}}\gamma(\epsilon)$. In particular, 
$a(\frac 1\epsilon)=\frac 1{\epsilon^{\frac 52}}\mathfrak q_{-}$ and $c_1=0$ if \eqref{erg-noi} holds.
 This proves Theorem \ref{main-tm}.
\end{proof}

The alternative splitting scheme \eqref{spl-sch1} can be treated by the same
strategy.  The estimates corresponding to \(Err_1^{k+1}\), \(Err_2^{k+1}\),
and the nonlinear part of \(Err_3^{k+1}\) are unchanged.  The only additional
term comes from replacing the exact stochastic convolution in \eqref{spl-sch}
by \(S_\lambda(\tau)\delta W_n\).  This term has zero first-order contribution,
and its weak error is controlled by a second-order Taylor expansion of the
Kolmogorov solution.  Therefore the proof needs the second derivative estimate
\eqref{time-ind2}; \revred{under the full non-degeneracy condition \eqref{erg-noi}, one
obtains the same type of bound, with the \(\mathfrak q_-\)-dependence appearing
only through the linear factor in \eqref{time-ind2}.}

\appendix

\section{Auxiliary proofs}
\label{app-auxiliary-proofs}

	\iffalse
	where $\eta$ is a sufficient small number, 
	\begin{align*}
		b_0&=\frac 1{\epsilon}(1-\eta)a_{2m+1} \frac {2} {2m+1}, \\
		b_{1,n}&=\left(\frac 1\epsilon + \lambda\right)C(\eta,m)\left(1+\|z_n\|^{2m+1}_E\right)+\frac 1{\epsilon} \frac{2m}  {m+q}(1-\eta)a_{2m+1}.
	\end{align*}
	Using iterations and taking $L^p(\Omega)$-norm yield that
	\begin{align*}
		\left\|v_{n+1}\right\|_{L^p(\Omega;E)}&\le   e^{-(c_0+b_0) (n+1)\tau}\|v_0\|_{L^p(\Omega;E)}+
		\sum_{k=0}^n e^{-(c_0+b_0) (k+1)\tau} \frac {1-e^{-b_0 \tau}}{b_0} \|b_{1,k}\|_{L^{p}(\Omega;E)}\\
		&\le e^{-(c_0+b_0) (n+1)\tau}\|v_0\|_{L^p(\Omega;E)}+
		\sum_{k=0}^n e^{-(c_0+b_0) (k+1)\tau} \tau \|b_{1,k}\|_{L^{p}(\Omega;E)}.
	\end{align*} 
	Applying the a priori bound of $z_n$ \eqref{pri-zn} and the discrete Gronwall's inequality, we complete the proof.
	\fi

\subsection{Proof of Lemma \ref{lem:negative-multiplier}}
\label{negative-multiplier}
\begin{proof}
We use the standard multiplier estimate: if \(s>d/2\) and \(0\le \gamma_1 \le s\), then
\[
  \|uv\|_{\HH^{\gamma_1}}
  \le C\|u\|_{\HH^{\gamma_1}}\|v\|_{H^s}.
\]
Hence, by duality, for \(g\in \HH^{\gamma_1}\) and \(y\in \HH^{-\gamma_1}\),
\[
\begin{aligned}
 \|gy\|_{\HH^{-s}}
 &= \sup_{\|\varphi\|_{\HH^s}\le1}
    |\langle y,g\varphi\rangle|  \\
 &\le \|y\|_{\HH^{-\gamma_1}}
    \sup_{\|\varphi\|_{\HH^s}\le1}\|g\varphi\|_{\HH^{\gamma_1}} \\
 &\le C\revpurple{\|g\|_{\HH^{\gamma_1}}}\|y\|_{\HH^{-\gamma_1}} .
\end{aligned}
\]
We apply this with
\[
  g(x)=-\epsilon^{-1}f'(x)+\lambda .
\]
Since \(f'\) is a polynomial of degree \(2m\), the Moser composition estimate gives
\[
 \|f'(x)\|_{\HH^{\gamma_1}}
 \le C\Bigl(1+\|x\|_E^{2m}
        +\revpurple{\|x\|_{\HH^{\gamma_1}}^{2m}}\Bigr).
\]
Moreover, because \(\lambda\) is fixed and \(0<\epsilon<1\),
\[
 \|\lambda y\|_{\HH^{-s}}
 \le C\lambda \|y\|_{\HH^{-s}}
 \le C\epsilon^{-1}\|y\|_{\HH^{-{\gamma_1}}} .
\]
Therefore,
\[
\begin{aligned}
\|(-\epsilon^{-1}f'(x)+\lambda)y\|_{\HH^{-s}}
&\le C\epsilon^{-1}\revpurple{\|f'(x)\|_{\HH^{\gamma_1}}}\|y\|_{\HH^{-\gamma_1}}
   +C\lambda\|y\|_{\HH^{-\gamma_1}} \\
&\le
C\epsilon^{-1}
\Bigl(1+\|x\|_E^{2m}
        +\|x\|_{\HH^{\gamma_1}}^{2m}\Bigr)
\|y\|_{\HH^{-\gamma_1}} .
\end{aligned}
\]
This proves the desired estimate.
\end{proof}

\subsection{Proof of  Lemma \ref{exa}}
\label{app-initial-free-E}

\begin{proof}

\textbf{Proof of (i):} We first prove \eqref{eq:initial-free-E0}.
Notice that $u$ can be decomposed into $y+z$,
	where 
	\begin{align}\label{random-eq}
		dy-(A-\lambda)y dt=f(y+z)+\lambda (y(t)+z(t))dt, \; y(0)=u(0).
	\end{align}
Here we denote $z$ the stochastic convolution, i.e., the mild solution of 
	$$dz-(A-\lambda)z dt=dW(t), \; z(0)=0.$$
	It follows from the similar arguments of the proof in \cite[Lemma 8.2.1 and Proposition 8.2.2]{Cer01} that for any $p\ge 1,$
	\begin{align}\label{sto-conv}
		\sup_{t\ge 0}\left\|z(t)\right\|_{L^p(\Omega;E)}+\sup_{t\ge 0}\left\|z(t)\right\|_{L^p(\Omega;\mathbb H^{\beta})}\le C_{p,\beta} \mathfrak q_\beta.
	\end{align}
	
By taking the subdifferential of the norm in $E$ on \eqref{random-eq}(see, e.g., \cite[Appendix]{Dap92}), and using the following  property of drift nonlinearity,
\begin{align}\label{con-theta0}
\frac 1\epsilon f(\xi)sign(\xi) \le -\frac {c}{\epsilon}|\xi|^{2m+1}+C\frac 1{\epsilon}
\end{align} 
for some postive constants $c\in (0,a_{2m+1})$ and $C>0,$
 we obtain that for $t\ge 0,$ $\zeta\in \partial \|y^{\delta}(t)\|_E$,
	
	\rev{\begin{align}\nonumber
		\frac {d^-}{dt}\left\|y(t)\right\|_E
		&\le \left\<(A-\lambda)y(t),\zeta\right\>+\left\<f(y+z)+\lambda y(t),\zeta\right\>\\\label{free-initial}
		&\le -\frac {c'}\epsilon \left\|y(t)\right\|_E^{2m+1}+ C\left(\lambda+\frac 1{\epsilon}\right)\left(1+\left\|z(t)\right\|_E^{2m+1}\right).
	\end{align}}
		
	Then Young's inequality and  Gronwall's inequality yield that 
	\begin{align*}
		\left\|y(t)\right\|_E\le  e^{-\frac {c_2}{\epsilon} t}\left\|u(0)\right\|_E+C(\lambda+\frac 1{\epsilon}) \int_{0}^t e^{-\frac {c_2}{\epsilon}(t-s)} \left(1+\left\|z(t)\right\|_E^{2m+1}\right)ds.
	\end{align*}
	
	Making use of \eqref{sto-conv}, we obtain that  for  $c_3>0$,
		\begin{align*}
			\left\|y\right\|_{L^{p}(\Omega; E)} \le e^{-\frac {c_3}{\epsilon} t} \left\|u(0)\right\|_E +C \mathfrak q_\beta^{2m+1}(1+\lambda \epsilon).
	\end{align*}
This, together with the fact that $\epsilon\in (0,1)$ and $\lambda$ is fixed, leads to \eqref{eq:initial-free-E0}.

Next we prove \eqref{eq:initial-free-E}. Fix $t>0.$ Similarly, taking the subdifferential on $\|y(t)\|_E^{p}$ for $p\ge 1$ and using again Young's inequality, we obtain 
\begin{align*}		\frac {d^-}{dt}\left\|y(t)\right\|_E^{p}
		&\le -\frac {c'}\epsilon \left\|y(t)\right\|_E^{2m+p}+ C\frac 1{\epsilon}\left(1+\left\|z(t)\right\|_E^{2m+1}\right)\|y(t)\|_E^{p-1}\\
		&\le -\frac {c_4}\epsilon \left\|y(t)\right\|_E^{2m+p}+C\frac 1{\epsilon}\left(1+\left\|z(t)\right\|_E^{(2m+1)p}\right).
	\end{align*}
	Taking expectation on the above inequality, we arrive at
\begin{align*}		\frac {d^-}{dt}\E \left\|y(t)\right\|_E^{p}
		&\le -\frac {c_4}\epsilon \E \left\|y(t)\right\|_E^{2m+p}+C\frac 1{\epsilon}\left(1+\E \left\|z(t)\right\|_E^{(2m+1)p}\right).
	\end{align*}		
Applying the standard scalar comparasion principle (similar as in the proof of \cite[Lemma 1.2.6]{Cer01}) to \eqref{free-initial} and using \eqref{sto-conv}, we obtain that 
\begin{align*}
 \E \left\|y(t)\right\|_E^{p}&\le C\mathfrak q_\beta^{\frac p{ {2m}+p}(2m+1)p}\left(1+\left(\frac {\epsilon}t\right)^{\frac p{2m}} \right).
\end{align*}
This, together with \eqref{sto-conv}, yields that 
\begin{align*}
  \left\|y(t)\right\|_{L^p(\Omega;E)}&\le C(\mathfrak q_\beta+\mathfrak q_\beta^{\frac {(2m+1)p}{2m+p}})\left(1+\left(\frac {\epsilon}t\right)^{\frac 1{2m}} \right),
\end{align*}
which compele the proof for \eqref{eq:initial-free-E} since $\frac {(2m+1)p}{2m+p}\le p.$

\textbf{Proof of (ii):}	We first prove \eqref{eq:initial-free-E2}. 
Since $\lambda>0$, the smoothing effect of $S_{\lambda}(t)$ \eqref{smooth0}-\eqref{smooth1} lead that there exists some $c_0\in (0,1),$ such that  for any $s\le t$,
any $\alpha<2$ and $\alpha_1\in [0,1],$
	\begin{align*}
		&\left\|y(t)-y(s)\right\|_{L^p(\Omega;\mathbb H)}\\
		&\le \Big\|S_{\lambda}(s)(S_{\lambda}(t-s)-I)y(0)\Big\|_{L^p(\Omega;\mathbb H)}\\
		&+ \Big\|\int_s^t S_{\lambda}(t-r) \left(f(y(r)+z(r))+\lambda (y(r)+z(r))\right)dr\Big\|_{L^p(\Omega;\mathbb H)}
		\\
		&+\Big\|\int_0^s S_{\lambda}(s-r) \Big(S_{\lambda}(t-s)-I\Big)\left( f(y(r)+z(r)+\lambda (y(r)+z(r)) \right) dr\Big\|_{L^p(\Omega;\mathbb H)}
		\\
		&\le C \|u(0)\|_{\mathbb H^{\beta}}|t-s|^{\frac \beta 2}
		\\
		&+C_{\alpha_1}\left(\frac {1}\epsilon+\lambda\right)|t-s|^{\alpha_1} \left(1+\sup_{r\in [s,t]}\left\|y(r)\right\|_{L^{p(2m+1)}(\Omega;E)}^{2m+1}+\sup_{r\in [s,t]}\left\|z(r)\right\|_{L^{p(2m+1)}(\Omega;E)}^{2m+1}\right)\\
		&+C_{\alpha}\left(\frac1\epsilon+ \lambda\right)\int_0^s e^{-c_0(s-r)}(s-r)^{-\frac \alpha 2} dr \Bigg(1+\sup_{r\in [s,t]}\left\|y(r)\right\|_{L^{p(2m+1)}(\Omega;E)}^{2m+1}\\
		&+\sup_{r\in [s,t]}\left\|z(r)\right\|_{L^{p(2m+1)}(\Omega;E)}^{2m+1}\Bigg) |t-s|^{\frac \alpha 2}\\
		&\le C_{\alpha,\alpha_1,\beta}\frac {1}\epsilon  (1+\mathfrak q_\beta^{(2m+1)^2}+\|u(0)\|_E^{2m+1}+\|u(0)\|_{\mathbb H^{\beta}})[ |t-s|^{\frac \alpha 2}+ |t-s|^{\alpha_1}+|t-s|^{\frac \beta 2}], 
	\end{align*}
	Combining this with $\alpha=\min(\beta,1), \alpha_1=\frac {\min(\beta,1)}2$, with the standard H\"older regularity estimate of the stochastic convolution (see, e.g., \cite[Chapter 6]{Cer01}) 
	\begin{align}\label{hold-convolution}
	\left\|z(t)-z(s)\right\|_{L^p(\Omega;\mathbb H)}\le C \mathfrak q_\beta |t-s|^{\min \left(\frac \beta 2,\frac 12\right)},
	\end{align} we can conclude that for $\beta<2$,
	\begin{align}\label{hold-udelta}
		\left\|u(t)-u(s)\right\|_{L^p(\Omega;\mathbb H)}\le C\mathfrak q_\beta^{(2m+1)^2} \frac {1}\epsilon (1+\|u(0)\|_E^{2m+1}+\|u(0)\|_{\mathbb H^{\beta}}) |t-s|^{\min \left(\frac \beta 2,\frac 12\right)}.
	\end{align}
This leads to \eqref{eq:initial-free-E2}.

Now we deal with \eqref{eq:initial-free-E1}.
	\iffalse
		\rev{\begin{align*}
		&\frac {d^-}{dt}\left\|y^{\delta}(t)\right\|_E\\
		&\le \left\<(A-\lambda)y^{\delta}(t),\zeta\right\>+\left\<\Theta_{\delta}^{\lambda}(y^{\delta}+z^{\delta})-\Theta_{\delta}^\lambda(z^{\delta}),\zeta\right\>+\left\<\Theta_{\delta}^\lambda(z^{\delta}),\zeta\right\>\\\
		&\le -\lambda \left\|y^{\delta}(t)\right\|_E+ C\left(\lambda+\frac 1{\epsilon}\right)\left(1+\left\|z^{\delta}(t)\right\|_E^{2m+1}\right)\\
		&+\int_0^1\Big\<\frac {-\frac {1}{\epsilon}g(\theta y^{\delta}(t)+z^{\delta}(t))+\lambda}{\left(1+\delta F(\theta y^{\delta}(t)+z^{\delta}(t))\right)^2}y^{\delta}(t),\zeta\Big\>d\theta\\
		&\le -\lambda \left\|y^{\delta}(t)\right\|_E+ C\left(\lambda+\frac 1{\epsilon}\right)\left(1+\left\|z^{\delta}(t)\right\|_E^{2m+1}\right).
	\end{align*}}
	\fi
	Similarly, we can obtain the spatial regularity estimate in the case $\beta<2$, that is for some $c_0\in (0,1),$
	\begin{align*}
		&\left\|y(t)\right\|_{L^p(\Omega;\mathbb H^{\beta})}\\
		&\le \left\|S_{\lambda}(t)u(0)\right\|_{L^p(\Omega;\mathbb H^{\beta})}+\Big\|\int_0^t S_{\lambda}(t-s)\left(f(y(s)+z(s))+\lambda(y(s)+z(s))\right)ds\Big\|_{\mathbb H^{\beta}}\\
		&\le C\left\|u(0)\right\|_{L^p(\Omega;\mathbb H^{\beta})}+C\left(\frac 1{\epsilon}+\lambda \right)\int_0^t e^{-c_0(t-s)}(t-s)^{-\frac \beta 2}ds\\
		& \quad\times \left(1+\sup_{s\in [0,t]}\left\|z(s)\right\|_{L^{p(2m+1)}(\Omega;E)}^{2m+1}+\sup_{s\in[0,t]}\left\|y(s)\right\|_{L^{p(2m+1)}(\Omega;E)}^{2m+1}\right)\\
		&\le C\left\|u(0)\right\|_{\mathbb H^{\beta}}+C\frac 1{\epsilon}\left(1+\left\|u(0)\right\|_{ E}^{2m+1}+\mathfrak q_\beta^{(2m+1)^2} \right),
	\end{align*}
	which, together with \eqref{sto-conv}, yields the spatial regularity estimate \eqref{eq:initial-free-E1}.

\end{proof}

\subsection{\texorpdfstring{Proof of \eqref{difference}:}{Proof of difference}}

\label{appendix-difference}

Let $\xi_0\in\mathbb R$ be fixed and let
\[
{Y(t)=\Phi_t^\lambda(\xi+\xi_0)-\xi_0.}
\]
Then $Y(0)=\xi$ and
\[
 \frac{d}{dt}Y(t)
 =-\frac1\epsilon f(Y(t)+\xi_0)+\lambda (Y(t)+\xi_0).
\]
By the dissipativity and polynomial growth of $f$, there exist constants
$c_0,C>0$, depending only on $m$ and $f$, such that for all $y,s\in\mathbb R$,
\begin{equation}\label{eq:scalar-diss}
\left(-\frac1\epsilon f(y+s)+\lambda(y+s)\right)\operatorname{sgn}(y)
\le
-\frac{c_0}{\epsilon}|y|+
C\frac1\epsilon
\left(1+|s|^{2m+1}\right).
\end{equation}
Consequently, in the sense of upper Dini derivatives,
\[
D^+|Y(t)|
\le
-\frac{c_0}{\epsilon}|Y(t)|
+
C\frac1\epsilon \left(1+|\xi_0|^{2m+1}\right).
\]
Solving this scalar differential inequality gives that for any $t>0,$
\begin{align}\label{eq:flow-shift}
|\Phi_t^\lambda(\xi+\xi_0)-\xi_0|
&\le
 e^{-\frac {c_0} \epsilon t}|\xi|  
+
C
\left(1+|\xi_0|^{2m+1}\right)(1-e^{-\frac {c_0}\epsilon t}),
\end{align}
which yields \eqref{difference}.

\subsection{Proof of Proposition \ref{reg-kl-short}}
\label{appendix-reg-finite}

\begin{proof}
The proof is similar to that of \cite[Proposition 4.1]{CH18}.
	Fix $T\ge t\ge s\ge 0$, define a two-parameter semigroup $II_{s,t}(\cdot)$ by
	\begin{align*}
		\revred{\partial_t II_{s,t} h= A II_{s,t} h-\frac 1\epsilon f' (u(t))\cdot II_{s,t} h, \; II_{s,s}h=h.}
	\end{align*}
	The energy estimate  immediately imply that 
	\begin{align}\label{ene-est}
		\frac {d}{dt}\left\|II_{s,t} h\right\|^2&\le 2\<A II_{s,t} h, II_{s,t} h\>+C\frac {1}{\epsilon} \left\|II_{s,t} h\right\|^2\\\nonumber 
		&\le -2\left\| \nabla II_{s,t} h\right\|^2+C\frac {1}{\epsilon} \left\|II_{s,t} h\right\|^2\\\nonumber 
		&\le C'\frac 1{\epsilon}\left\| II_{s,t} h\right\|^2.
	\end{align}
	This yields that 
	\begin{align}\label{pri-reg1}
		\|II_{s,t} h\|^2\le \exp\Bigg(C'\frac 1{\epsilon}(t-s)\Bigg) \|h\|^2, \; s\le t\le T.
	\end{align}
	Recall \eqref{var-1}.
	Applying the decomposition $\widetilde \eta^h(t,x)=\eta^h(t,x)-S_{\lambda}(t)h,$
	it follows that 
	\begin{align*}
		\revred{\frac {d \widetilde \eta^h}{dt}=(A-\lambda)\widetilde \eta^h+\left[-\frac 1\epsilon f'\left(u(t,x)\right)+\lambda\right]\cdot \widetilde \eta^h+\left[-\frac 1\epsilon f'\left(u(t,x)\right)+\lambda\right]\cdot S_\lambda(t) h,}
	\end{align*}
	where $\widetilde \eta^h(0)=0$.
	Notice that 
	\begin{align*}
		\revred{\widetilde \eta^h(t)}&\revred{=\int_0^t II_{s,t} \left(\left[-\frac 1\epsilon f'\left(u(s,x)\right)+\lambda\right]\cdot S_{\lambda}(s)h\right)ds.}
	\end{align*}
	Thanks to  the smoothing effect \eqref{smooth0} of  $S_{\lambda}(\cdot)$,  by \eqref{pri-reg1} and \eqref{eq:initial-free-E0}, we have that for some $c_0\in (0,\lambda)$, and any $\alpha_1\in (0,1),$
	\begin{align*}
		&\left\|\widetilde \eta^h(t)\right\|_{\mathbb H}\revred{\le} \int_0^t \left\| II_{s,t} \left(\revred{\left[-\frac 1\epsilon f'\left(u(s,x)\right)+\lambda\right]} \cdot S_{\lambda}(s)h\right)\right\| ds\\
		&\le C \int_{0}^t \exp\left(C'\frac 1{\epsilon}(t-s) \right)\left\|\revred{\left[-\frac 1\epsilon f'\left(u(s,x)\right)+\lambda\right]} \right\|_E  e^{-c_0  s} s^{-\alpha_1}ds \|h\|_{\mathbb H^{-2\alpha_1}}\\
		&\le  C \frac1\epsilon\exp\left(C'\frac 1\epsilon t\right)P_{2m+1}(\mathfrak q_{\beta})^{2m} \left(1+\|x\|_E^{2m}\right)  \|h\|_{\mathbb H^{-2\alpha_1}},
	\end{align*}
	where we use the fact that $\int_0^{\infty}e^{-c_0  s} s^{-\alpha_1}ds<\infty.$
	
	This, together with \eqref{smooth0}, we conclude that 
	\begin{align*}
		\left\|\eta^h(t,x)\right\|_{L^p(\Omega;\mathbb H)}&\le   C \frac1\epsilon\exp\left(C'\frac 1\epsilon t\right) P_{2m+1}(\mathfrak q_{\beta})^{2m}  \left(1+\|x\|_E^{2m}\right)  \|h\|_{\mathbb H^{-2\alpha_1}}\\
		&\quad \times 
		\left(1+e^{-c_0 t}t^{-\alpha_1}\right).
	\end{align*}
	which yields \eqref{reg-kol-1} thanks to \eqref{rep-kol}.
	
	This also implies that 
	\begin{align}\nonumber
		\|II_{s,t} h \|_{L^p(\Omega;\mathbb H)} &\le C\frac1\epsilon\exp\left(C'\frac 1\epsilon (t-s)\right)P_{2m+1}(\mathfrak q_{\beta})^{2m}   \left(1+\|x\|_E^{2m}\right)  \|h\|_{\mathbb H^{-2\alpha_1}}\\\label{smo-two}
		&\quad \times
		\left(1+e^{-c_0(t-s)}(t-s)^{-\alpha_1}\right).
	\end{align}

	Similarly,  recalling \eqref{var-2}, according to the fact that 
	\begin{align}\label{4.9}
		\revred{\zeta^{h_1,h_2}(t,x)=-\int_0^t II_{s,t} \left(\frac 1\epsilon f''(u(s,x))\cdot \left(\eta^{h_1}(s,x),\eta^{h_2}(s,x)\right)\right)ds,}
	\end{align}
	and \eqref{smo-two}, applying the Sobolev embedding theorem $\mathbb H^{2\alpha}\hookrightarrow L^{\infty}$ for $\alpha\in (\frac d4,1),$ and H\"older's inequality,
	it follows that 	\begin{align*}
		&\left\|\zeta^{h_1,h_2}(t,x)\right\|_{L^p(\Omega;\mathbb H)}\\
		&\le C \frac 1\epsilon\int_0^t \exp\left(C'\frac {1}{\epsilon}(t-s) \right)P_{2m+1}(\mathfrak q_{\beta})^{2m}   \left(1+\|x\|_E^{2m}\right) \left(1+e^{-c_0(t-s)}(t-s)^{-\alpha}\right) \\\nonumber 
		&\quad 
		 \Big\|f''(u(s,x))\cdot \big(\eta^{h_1}(s,x),  \eta^{h_2}(s,x)\big)\Big\|_{L^{2p}(\Omega;\mathbb H^{-2\alpha})} ds\\
		&\le C (\frac 1\epsilon)^2\int_0^t \exp\left(C'\frac {1}{\epsilon}(t-s) \right) \left(1+\|x\|_E^{4m-1}\right)  \left(1+e^{-c_0(t-s)}(t-s)^{-\alpha}\right) \\\nonumber 
		&\quad 
		P_{2m+1}(\mathfrak q_{\beta})^{4m-1}   \left\|\eta^{h_1}(s,x)\right\|_{L^{6p}(\Omega;\mathbb H)} \left\| \eta^{h_2}(s,x)\right\|_{L^{6p}(\Omega;\mathbb H)}  ds.
	\end{align*}
	Therefore, applying \eqref{smo-two}  again and \eqref{exa}, we obtain that 
	\begin{align*}
		&\left\|\zeta^{h_1,h_2}(t,x)\right\|_{L^{p}(\Omega;{\mathbb H})}\\
		&\le C\left(\frac1\epsilon\right)^{4}
		\exp\left(\frac {C't}{\epsilon} \right)\int_0^t  \left(1+\|x\|_E^{8m-1}\right) \left(1+e^{-c_0(t-s)}(t-s)^{-\alpha}\right)\\
		&
		\quad \times P_{2m+1}(\mathfrak q_{\beta})^{8m-1}  \left(1+e^{-c_0s}s^{-\alpha_1}\right) \left(1+e^{-c_0 s}s^{-\alpha_2}\right)  ds \|h_1\|_{\mathbb H^{-2\alpha_1}} \|h_2\|_{\mathbb H^{-2\alpha_2}}.
	\end{align*}
	From the fact that $\int_{0}^1 \xi^{-\alpha}(1-\xi)^{-\alpha_1-\alpha_2} d\xi<\infty$ with $\alpha_1+\alpha_2\in(0,1),\alpha\in (\frac d4,1),$  it follows that
	\begin{align}\nonumber
		&\left\|\zeta^{h_1,h_2}(t,x)\right\|_{L^{p}(\Omega;{\mathbb H})}\\\nonumber
		&\le C  P_{2m+1}(\mathfrak q_{\beta})^{8m-1}  \left(\frac1\epsilon\right)^{4}
		\exp\left(\frac {C't}{\epsilon} \right)\left(1+\|x\|_E^{8m-1}\right)(1+t^{1-\alpha-\alpha_1-\alpha_2})\\\label{var-sec}
		&\quad \times \|h_1\|_{\mathbb H^{-2\alpha_1}}\|h_2\|_{\mathbb H^{-2\alpha_2}}.
	\end{align}
	
	By combining the above estimates and \eqref{rep-kol1} together, we complete the proof of \eqref{reg-kol-2} due to $\phi\in\mathcal C_b^2(\mathbb H).$
\end{proof}

\subsection{\texorpdfstring{Proof of \eqref{asy-strong-feller} in Lemma \ref{Strong Feller property} under the condition that $\gamma(\epsilon)<+\infty$}{Proof of the asymptotic strong Feller estimate under the low-mode non-degeneracy condition}}
\label{general-time-ind}

\begin{proof}
We now prove the estimate \eqref{asy-strong-feller} under the low-mode non-degeneracy condition
\(\gamma(\epsilon)<\infty\), following the asymptotic coupling argument of
\cite[Lemma 4.1]{CW2024}.

\textbf{Step 1: Exponential estimate for asymptotic coupling}

Let $u(t)=u(t,x)$ be the solution of \eqref{sac} starting from $x$. For another initial value $y$, define a coupled process $\widetilde u(t)$ by
\begin{equation}\label{eq:coupled-lowmode}
    d\widetilde u(t)
    +\left(
        A\widetilde u(t)+\epsilon^{-1}f(\widetilde u(t))
        +\frac{1+L_f}{\epsilon}\pie(\widetilde u(t)-u(t))
      \right)d t
    =d W(t),
    \qquad \widetilde u(0)=y.
\end{equation}
Set
\[
    \widehat Y(t):=\widetilde u(t)-u(t).
\]
Then
\begin{equation}\label{eq:Z-eq}
d\widehat  Y(t)
    +A\widehat  Y(t)dt+\epsilon^{-1}\bigl(f(\widetilde u(t))-f(u(t))\bigr)dt
    +\frac{1+L_f}{\epsilon}\pie \widehat Y(t)dt=0.
\end{equation}
Taking the inner product with $\widehat Y(t)$ and using the one-sided Lipschitz continuity of $f$, we obtain
\begin{align}\label{eq:Z-energy-1}
    \frac12 d\norm{\widehat Y(t)}^2
    &
    \le
      -\norm{\nabla \widehat Y(t)}^2 dt
      +\frac{L_f}{\epsilon}\norm{ \widehat Y(t)}^2dt 
      -\frac{1+L_f}{\epsilon}\norm{\pie  \widehat Y(t)}^2dt.
\end{align}
By the definition of $N_\epsilon$,
\begin{equation}\label{eq:spectral-gap}
    \norm{\widehat Y(t)}_{\mathbb H^{1}}^2
    \ge \lambda_{N_\epsilon+1}\norm{(I-\pie)\widehat Y(t)}^2
    \ge \frac{1+L_f}{\epsilon}\norm{(I-\pie)\widehat Y(t)}^2.
\end{equation}
Combining \eqref{eq:Z-energy-1} with \eqref{eq:spectral-gap} gives
\begin{align*}
   & \frac12d\norm{\widehat Y(t)}^2\\
    &\le
    -\frac{1+L_f}{\epsilon}\norm{(I-\pie)\widehat Y(t)}^2dt
    +\frac{L_f}{\epsilon}\norm{\widehat Y(t)}^2dt
    -\frac{1+L_f}{\epsilon}\norm{\pie \widehat Y(t)}^2 dt\\
    &=
    -\frac1\epsilon\norm{(I-\pie)\widehat Y(t)}^2dt
    -\frac1\epsilon\norm{\pie \widehat Y(t)}^2dt
    =-\frac1\epsilon\norm{\widehat Y(t)}^2dt.
\end{align*}
Therefore
\begin{equation}\label{eq:Z-decay}
    \norm{\widetilde u(t)-u(t)}
    \le e^{-t/\epsilon}\norm{x-y},
    \qquad t\ge0.
\end{equation}
This is the source of the asymptotic term in \eqref{asy-strong-feller}. 
%Notice that this is not a terminal coupling: the two solutions need not coincide at time $t$.

%\medskip
\textbf{Step 2: Girsanov transform on the low modes.}
Rewrite \eqref{eq:coupled-lowmode} as the original equation driven by a shifted Wiener process. Define
\begin{equation}\label{eq:xi-def}
    \xi_s:=\frac{1+L_f}{\epsilon}(Q^{1/2})^{\dagger}\pie(\widetilde u(s)-u(s)).
\end{equation}
where $(Q^{1/2})^{\dagger}$ is the pseudo-inverse on the span of those $e_i$ with $q_i>0$. Since $\pie$ only involves $1\le i\le N_\epsilon$, the quantity \eqref{eq:xi-def} is well-defined. By \eqref{eq:Neps} and \eqref{eq:Z-decay},
\begin{equation}\label{eq:xi-bound}
    \norm{\xi_s}^2
    \le
    \gamma(\epsilon)^2\left(\frac{1+L_f}{\epsilon}\right)^2
    e^{-2s/\epsilon}\norm{x-y}^2.
\end{equation}
Hence
\begin{equation}\label{eq:xi-integral}
    \int_0^t\norm{\xi_s}^2d s
    \le
    \frac{\gamma(\epsilon)^2(1+L_f)^2}{2\epsilon}
\norm{x-y}^2.
\end{equation}
Let
\begin{equation}\label{eq:R-density}
    R_t
    :=\exp\left(
        \int_0^t\ip{\xi_s}{d W_s}
        -\frac12\int_0^t\norm{\xi_s}^2d s
      \right).
\end{equation}
By Girsanov's theorem, under the probability measure $R_t\mathbb P$, the process
\[
    \widetilde W(t):=W(t)-\int_0^t Q^{1/2}\xi_s ds
\]
is again a $Q$-Wiener process. Consequently, by weak uniqueness, $\widetilde u(t)$ under $R_t\mathbb P$ has the same law as the original solution starting from $y$.

\medskip
\noindent\textbf{Step 3: Asymptotic log-Harnack inequality.}
Let $\psi>0$ be bounded and Lipschitz. Then
\begin{align}\label{eq:logharnack-start}
    P_t\log\psi(y)
    &=\E\bigl[R_t\log\psi(\widetilde u(t))\bigr] \notag\\
    &\le
      \E\bigl[R_t\log\psi(u(t))\bigr]
      +\norm{D\log\psi}_{\infty}\E\bigl[R_t\norm{\widetilde u(t)-u(t)}\bigr].
\end{align}
Using \eqref{eq:Z-decay}, the second term is bounded by
\[
    \norm{D\log\psi}_{\infty}e^{-t/\epsilon}\norm{x-y}.
\]
For the first term, the entropy inequality gives
\begin{equation}\label{eq:entropy}
    \E\bigl[R_t\log\psi(u(t))\bigr]
    \le
    \log P_t\psi(x)+\E[R_t\log R_t].
\end{equation}
Moreover, under $R_t\mathbb P$,
\begin{equation}\label{eq:entropy-cost}
    \E[R_t\log R_t]
    =\frac12\E_{R_t\mathbb P}\int_0^t\norm{\xi_s}^2d s
    \le
    \frac{\gamma(\epsilon)^2(1+L_f)^2}{4\epsilon}\norm{x-y}^2.
\end{equation}
Combining \eqref{eq:logharnack-start}--\eqref{eq:entropy-cost}, we obtain
\begin{equation}\label{eq:asymp-logharnack}
    P_t\log\psi(y)
    \le
    \log P_t\psi(x)
    +\norm{D\log\psi}_{\infty}e^{-t/\epsilon}\norm{x-y}
    +\frac{\gamma(\epsilon)^2(1+L_f)^2}{4\epsilon}\norm{x-y}^2.
\end{equation}
This is the asymptotic log-Harnack inequality for low-mode coupling.

\medskip
\noindent\textbf{Step 4: Gradient estimate.}
The standard implication from the asymptotic log-Harnack inequality \eqref{eq:asymp-logharnack} \cite{MR4013873} gives, for every bounded Lipschitz $\phi$,
\begin{equation}\label{eq:gradient-final}
    \norm{DP_t\phi}_{\infty}
    \le
    e^{-t/\epsilon}\norm{D\phi}_{\infty}
    +\frac{(1+L_f)\gamma(\epsilon)}{2\sqrt{\epsilon}}\norm{\phi}_{\infty}.
\end{equation}
This proves \eqref{asy-strong-feller} with $c_0=c_1=1,$ and 
\begin{align}\label{bound-case1}
b(\epsilon^{-1},t)=\frac{(1+L_f)\gamma(\epsilon)}{2\sqrt{\epsilon}}.
\end{align}

Next, we extend this estimate to measurable $\phi$ of polynomial growth and $\phi\in \mathcal C_b^1(\mathbb H)\cap B_{V,q}(\mathbb H)$. 
By the Markov property, it holds that $P_t=P_{\frac t2}P_{\frac t2}\phi.$
By \eqref{eq:gradient-final} and \eqref{eq:initial-free-E}, we have that 
\begin{align*}
    \norm{DP_t\phi}_{\infty}
    &\le
    e^{-\frac {t}{2\epsilon}}\norm{DP_{\frac  t2}\phi}_{\infty}
    +\frac{(1+L_f)\gamma(\epsilon)}{2\sqrt{\epsilon}}\norm{P_{\frac 12 t}\phi}_{\infty}\\
    &\le e^{-\frac {t}{2\epsilon}}\norm{DP_{\frac  t2}\phi}_{\infty}
    +\frac{(1+L_f)\gamma(\epsilon)}{2\sqrt{\epsilon}}P_{q}(\mathfrak q_{\beta})^q\left(1+\left(\frac{\epsilon}{t}\right)^{q/(2m)}\right)\\
    &\le  e^{-\frac {t}{\epsilon}}\norm{D\phi}_{\infty}+\frac{(1+L_f)\gamma(\epsilon)}{\sqrt{\epsilon}}P_{q}(\mathfrak q_{\beta})^q \left(1+\left(\frac{\epsilon}{t}\right)^{q/(2m)}\right).
\end{align*}
This leads to \eqref{asy-strong-feller} with $c_0=c_1=1,$ 
and 
\begin{align}\label{polynomial-asy}
\tilde b(\epsilon^{-1},t,\mathfrak q_{\beta},\gamma(\epsilon))
&=
\frac{(1+L_f)\gamma(\epsilon)}{\sqrt{\epsilon}}
P_{q}(\mathfrak q_{\beta})^q
\left(1+\left(\frac{\epsilon}{t}\right)^{q/(2m)}\right).
\end{align}

\end{proof}

\bibliographystyle{plainnat}
\bibliography{bib}

\end{document}